\documentclass[10pt]{article}

\usepackage{graphicx,hyperref}
\usepackage{amsmath}
\usepackage[version=4]{mhchem}
\usepackage{siunitx}
\usepackage{longtable,tabularx}
\usepackage{pstricks}
\usepackage{amsfonts}
\usepackage{latexsym}
\usepackage{epsfig}
\usepackage{subfig}
\usepackage{framed,fancybox}
\usepackage{bbm}
\usepackage{multirow}
\usepackage[T1]{fontenc}
\usepackage{threeparttable}
\usepackage{rotating}
\usepackage{amssymb,amscd}
\makeatletter
\def\myVCENTER#1{\vcenter{\hbox{$\m@th#1$}}}
\makeatother

\long\def\symbolfootnote[#1]#2{\begingroup\def\thefootnote{\fnsymbol{footnote}}\footnote[#1]{#2}\endgroup}

\definecolor{shadecolor}{gray}{0.99}
{\endMakeFramed}

\definecolor{shadecolor}{gray}{0.99}
\long\def\symbolfootnote[#1]#2{\begingroup\def\thefootnote{\fnsymbol{footnote}}\footnote[#1]{#2}\endgroup}
\def\qed{\hfill{$\vcenter{\hrule height1pt \hbox{\vrule width1pt height5pt
    \kern5pt \vrule width1pt} \hrule height1pt}$} \medskip}
\newcommand{\m}[1]{{\bf{#1}}}
\newcommand{\g}[1]{\boldsymbol #1}
\newcommand{\bb}[1]{\mathbb #1}
\newcommand{\tr}{^{\sf T}}
\newcommand{\C}[1]{{\cal {#1}}}

\title{\bf Comparison of Derivative Estimation Methods in Solving Optimal Control Problems Using Direct Collocation}

\author{Yunus M.~Agamawi\thanks{Ph.D.~Student, Department of Mechanical and Aerospace Engineering.  E-mail:  yagamawi@ufl.edu.} \\ Anil V.~Rao\thanks{Associate Professor, Department of Mechanical and Aerospace Engineering, Erich Farber Faculty Fellow and University Term Professor.  E-mail: anilvrao@ufl.edu.  Associate Fellow, AIAA.  Corresponding Author.} \vspace{12pt} \\ {\em University of  Florida} \\ {\em Gainesville, FL, 32611-6250}}

\date{}
\begin{document}

\renewcommand{\baselinestretch}{1}
\normalsize\normalfont 
\maketitle
\begin{abstract}
  \noindent A study is conducted to evaluate four derivative estimation methods when solving a large sparse nonlinear programming problem that arises from the approximation of an optimal control problem using a direct collocation method.  In particular, the Taylor series-based finite-difference, bicomplex-step, and hyper-dual derivative estimation methods are evaluated and compared alongside a well known automatic differentiation method.  The performance of each derivative estimation method is assessed based on the number of iterations, the computation time per iteration, and the total computation time required to solve the nonlinear programming problem.  The efficiency of each of the four derivative estimation methods is compared by solving three benchmark optimal control problems.  It is found that while central finite-differencing is typically more efficient per iteration than either the hyper-dual or bicomplex-step, the latter two methods have significantly lower overall computation times due to the fact that fewer iterations are required by the nonlinear programming problem when compared with central finite-differencing.  Furthermore, while the bicomplex-step and hyper-dual methods are similar in performance, the hyper-dual method is significantly easier to implement.  Moreover, the automatic differentiation method is found to be substantially less computationally efficient than any of the three Taylor series-based methods.  The results of this study show that the hyper-dual method offers several benefits over the other three methods both in terms of computational efficiency and ease of implementation.
\end{abstract}

\section*{Nomenclature}


{\renewcommand\arraystretch{1.0}
\noindent\begin{longtable*}{lcl}
$\m{a}$	& = & vector field for right-hand side of dynamics \\
$\m{A}$	& = &  matrix defining vector field for right-hand side of dynamics at collocation points \\
$AD$		& = & automatic differentiation method	\\
$\m{b}$	& = & vector field for event constraints \\
$BC$		& = & bicomplex-step method	\\
$\m{c}$	& = & vector field for path constraints \\
$\m{C}$	& = &  matrix defining vector field for path constraints at collocation points \\
$\m{D}$	& = &  Legendre-Gauss-Radau differentiation matrix \\
$\m{e}$	& = & endpoint vector \\
$\m{E}$	& = &  endpoint approximation vector \\
$EC$		& = & exact sparsity central finite-difference method	\\
$\m{g}$	& = & vector field for integrands \\
$\m{G}$	& = &  matrix defining vector field for integrands at collocation points \\
$h$			& = & perturbation step size	\\
$h_{\mathcal{O}_1}$ & = &  base perturbation step size for finite-difference	\\
$HD$		& = & hyper-dual method	\\
$I$			& = & number of NLP solver iterations	\\
$\C{J}$	& = & objective functional \\
$K$		& = & number of mesh intervals used in phase	\\
$\ell_j^{(k)}$ & = & Lagrange polynomial $j$ of mesh interval $k$ \\
$n_y$	& = & number of state components	\\
$n_u$	& = & number of control components 	\\
$n_q$	& = & number of integral components 	\\
$n_c$	& = & number of path constraints	\\
$n_b$	& = & number of event constraints	\\
$n_s$	& = & number of static parameters	\\
$N_k$	& = & number of collocation points used in mesh interval $k$	\\
$N$		& = & total number of collocation points used in phase	\\
$OC$		& = & over-estimated sparsity central finite-difference method	\\
$P$		& = & number of phases in problem	\\
$\m{q}$	& = & integral vector \\
$\m{Q}$	& = &  integral approximation vector \\
$\m{s}$	& = & static parameters vector	\\
$\C{S}_k$ & = & mesh interval $k$\\
$t_0$	& = & initial time \\
$t_f$ 	& = & final time \\
$\C{T}$		& = & computation time	\\
$\m{u}(\tau)$	& = & control \\
$\m{U}$	& = &  matrix of control parameterization at collocation points \\
$\m{y}(\tau)$	& = & state \\
$\m{Y}(\tau)$ & = & state approximation	\\
$\m{Y}$	& = &  matrix of state approximation at discretized points \\
$\Gamma_{\xi}^{\mathcal{A}}$ & = & percent change of criteria $\xi\in\{I,\C{T},\Phi\}$ for method \\
 & & $\mathcal{A}\in\{EC,BC,HD,AD\}$ relative to $OC$ method \\
$\g{\Delta}$	& = &  defect constraints matrix \\
$\g{\rho}$	& = & integral approximation constraints vector \\
$\tau_j^{(k)}$ & = & support point $j$ of mesh interval $k$ \\
$\Phi$		& = & average time per NLP iteration \\
\end{longtable*}}

\renewcommand{\baselinestretch}{1}
\normalsize\normalfont 

\section{Introduction \label{sect:intro}}

Optimal control problems arise in a wide variety of engineering and non-engineering disciplines.  Due to the increasing complexity of modern optimal control problems and the inability to solve these problems analytically, in the past few decades the focus of research in optimal control has transitioned from theory to computation.  One approach that has become important in the numerical solution of optimal control problems is the class of {\em direct collocation methods}.  In a direct collocation method, the state and control of the optimal control problem are approximated using a set of basis or trial functions and the constraints are enforced at a special set of points called nodes.  The continuous optimal control problem is then approximated as a large sparse nonlinear programming problem (NLP) \cite{Betts3} and this sparse NLP is solved using well known software \cite{Gill3,Biegler2}.

The sparse NLP arising from a direct collocation method is generally solved using a gradient-based method.  The two categories of gradient-based methods for solving NLPs are quasi-Newton methods and Newton methods.  In a quasi-Newton gradient-based method, the NLP objective function gradient and constraints Jacobian are required along with a quasi-Newton approximation to the inverse of the Hessian of the NLP Lagrangian (for example, BFGS or DFP quasi-Newton approximations).  A Newton method adds to a quasi-Newton method the requirement that the Lagrangian Hessian be computed.  Perhaps the most well known NLP solver that implements a quasi-Newton method is {\em SNOPT} \cite{Gill1}, while perhaps the most well known Newton NLP solver is {\em IPOPT} \cite{Biegler2}.   In order to make it tractable to solve the NLP using a gradient-based method, it is necessary to efficiently and accurately compute the required first- and second-order derivatives of the NLP functions.  In fact, the vast majority of time required for the NLP solver to obtain a solution is spent computing the required derivatives.  It is also important to note that, even if the derivatives are computed efficiently, inaccuracies in the derivatives can lead to poor convergence or non-convergence of the NLP solver.  Thus, the ability to supply efficient and accurate derivatives is of utmost importance when solving an NLP arising from the approximation of an optimal control problem using direct collocation.  

A variety of different approaches have been developed for estimating derivatives.  These approaches can be divided broadly into algorithmic methods and approximation methods.  In an algorithmic method, a derivative is obtained by decomposing the function into a sequence of elementary operations and applying rules of differential calculus.  Algorithmic approaches include analytic differentiation, symbolic differentiation, and automatic differentiation \cite{Martins2}.  Analytic differentiation is the process of obtaining the derivative manually by differentiating each step in the chain rule by hand \cite{Sobieszczanski1}.  As a result, analytic differentiation is often referred to as hand differentiation.  Symbolic differentiation is the process of obtaining a derivative using computer algebra such as that implemented in {\em Maple} or {\em Mathematica}.  Finally, automatic differentiation is the process of decomposing a computer program into a sequence of elementary function operations and applying the calculus chain rule algorithmically through the computer \cite{Griewank0}.  It is noted that all of the aforementioned algorithmic methods can be derived from the unifying chain rule described in Ref.~\cite{Martins2}.  Finally, an approximation method evaluates the function at one or more neighboring points and estimates the derivative using this information.

Algorithmic and approximation methods for computing derivatives have relative advantages and disadvantages.  The key advantage of an algorithmic method is that an exact derivative (that is, a derivative accurate to machine precision on a computer) is obtained.  In addition, in some cases using an algorithmic method makes it possible to obtain a derivative in a computationally efficient manner, thus maximizing the reliability and computational efficiency of solving the nonlinear programming problem arising from a direct collocation method.  It is noted, however, that algorithmic methods have several limitations.  First, for a practical problem it is generally intractable to obtain an analytic derivative.  Second, symbolic differentiation used in computer algebra systems suffer from expression explosion that make it intractable to compute derivatives whose forms are algebraically complex.  Third, automatic differentiation can only be performed on functions where the differentiation rule has been defined via a library or database of functions.  In fact, it is often intractable to define derivative rules for highly complicated functions (for example, determining the derivative rule associated with the interpolation of a multi-dimensional table).

In particular, the well-known open source automatic differentiation software ADOL-C \cite{Griewank1,Walther1}, which is included in this paper for comparison as the automatic differentiation method, requires a tape recording of the function evaluation for each derivative calculation.  Because the record used in ADOL-C does not account for changes in flow control, problems utilizing piecewise models or interpolation of tabular data must be retaped upon each input in order to ensure the proper derivative calculation is performed.  Although taping the record of the function evaluation for the derivative to be calculated can be efficient if reused, the conditional flow control nature of many optimal control problems requires retaping in order to ensure the proper flow control is computed for each input.  A key advantage of an approximation method is that only the function needs to be evaluated in order to estimate the derivative which generally makes an approximation method computationally efficient.  It is noted, however, that typical approximation methods (for example, finite-differencing) are not very accurate, thereby reducing the robustness with which a solution to the NLP can be obtained.

While a general NLP may not contain a great deal of structure, it has been shown in Refs.~\cite{Betts13,Patterson2012,Agamawi2018} that the derivative functions of an NLP arising from the direct collocation of an optimal control problem have a very well-defined structure.  In particular, Refs.~\cite{Betts13,Patterson2012,Agamawi2018} have shown that the derivative functions of the NLP arising from a direct collocation method can be obtained with increased efficiency by computing the derivatives of the functions associated with the optimal control problem at the collocation points as opposed to computing the derivatives of the NLP functions.  This exploitation of sparsity in a direct collocation NLP makes it possible to compute the derivatives of the lower-dimensional optimal control problem functions relative to the much higher-dimensional derivatives of the NLP functions.  Given the increased efficiency that arises from sparsity exploitation, any derivative method employed when solving an optimal control problem using direct collocation should exploit this NLP sparsity.  

The aforementioned discussion makes it relevant and useful to evaluate the effectiveness of using various approximation methods for computing derivatives required when solving an NLP arising from a direct collocation method.  Typically, finite-difference approximation methods are used when solving a direct collocation NLP~\cite{Betts3,Gill3,Betts2}.  In recent years, however, much more accurate derivative approximation methods have been developed.  Two such recent developments that enable the computation of higher-order (that is, beyond first-derivative) derivatives are the bicomplex-step~\cite{Lantoine1}, and hyper-dual~\cite{Fike2011} derivative approximations.  In a bicomplex-step derivative approximation, a new class of numbers, called {\em bicomplex numbers}, is created.  As its name implies, a bicomplex number is one that has one real component, two imaginary components, and one bi-imaginary component.  Derivative approximations are then obtained by evaluating the functions at a particular bicomplex input argument.  While more accurate than finite-difference methods, the bicomplex-step derivative approximation still has inaccuracies stemming from the bicomplex arithmetic \cite{Lantoine1}.  Next, in a hyper-dual derivative approximation, a new class of numbers, called {\em hyper-dual numbers}, is created \cite{Fike2011}.  As its name implies, a hyper-dual number is one that has one real component, two imaginary components, and one bi-imaginary component.  While hyper-dual numbers share some properties with those of bicomplex numbers, the mathematics associated with a hyper-dual number are fundamentally different from those of a bicomplex number.  The difference in the mathematics associated with a hyper-dual number relative to other approximation methods leads to a derivative approximation that does not suffer from truncation error and is insensitive to the step size.  Finally, it is noted that all three of the aforementioned derivative approximation methods fall into the category of Taylor series-based approximation methods.

The objective of this research is to evaluate the effectiveness of using the following four derivative estimation methods for use when solving an NLP arising from the direct collocation of an optimal control problem: central finite-differencing, bicomplex-step, hyper-dual, and automatic differentiation.  Specifically, the objective of this paper is to provide a comprehensive comparison of the computational efficiency with which a direct collocation NLP can be solved using a Newton NLP solver with the four aforementioned derivative estimation methods.  The evaluation is based on the number of iterations required to obtain a solution to the NLP, the overall computation time required to solve the NLP, and the average computation time required per iteration to compute the required derivatives.  It is found that while central finite-differencing is typically more efficient per iteration than either the hyper-dual or bicomplex-step, the latter two methods have a significantly lower overall computation time due to the fact that fewer iterations are required by the nonlinear programming problem solver using these latter two methods when compared with central finite-differencing.  Furthermore, it is found that although the bicomplex-step and hyper-dual are similar in performance, the hyper-dual derivative approximation is significantly easier to implement.  Moreover, the three Taylor series-based derivative approximations are found to be substantially more computationally efficient methods when compared to an automatic differentiation method provided by the open source software ADOL-C \cite{Griewank1,Walther1}.  The results of this study show that the hyper-dual derivative approximation offers several benefits over the other two derivative approximations and the tested automatic differentiation derivative estimation method, both in terms of computational efficiency and ease of implementation.

The contributions of this research are as follows.  First, this research provides a systematic evaluation of the performance of using various Taylor series-based derivative approximation methods and an automatic differentiation method when solving a nonlinear programming problem arising from a direct collocation method for optimal control.  In addition, a method for determining the exact sparsity patterns of the first- and second-order NLP derivative functions is developed by estimating the continuous-time optimal control problem functions using bicomplex-step or hyper-dual derivative approximations.  Furthermore, the implementation and performance behavior of the NLP solver when utilizing the three aforementioned Taylor series-based derivative approximation methods and automatic differentiation as methods are compared against one another.  As mentioned, the results of this paper show that, among the methods evaluated, the hyper-dual derivative approximation is most effective and has the greatest ease of implementation for solving optimal control problems using direct collocation methods.

This paper is organized as follows.  Section~\ref{sect:multi-phase} defines the general multiple-phase optimal control problem.  Section~\ref{sect:LGR} describes the rationale for using Legendre-Gauss-Radau collocation points as the set of nodes to discretize the continuous optimal control problem and the Legendre-Gauss-Radau collocation method.  Section~\ref{sect:NLP} briefly overviews the form of the nonlinear programming problem that results from transcribing the optimal control problem.  Section~\ref{sect:deriv-approx} describes the basis of the four derivative estimation methods: central finite-differencing, bicomplex-step derivative approximation, hyper-dual derivative approximation, and automatic differentiation.  Section~\ref{sect:DerivativeComparison} briefly compares the derivative approximations obtained when using central finite-differencing versus a bicomplex-step derivative approximation, hyper-dual derivative approximation, or automatic differentiation.  Section~\ref{sect:Examples} provides examples that demonstrate the performance of the NLP solver when using a central finite-difference method versus a bicomplex-step method, hyper-dual method, or automatic differentiation method.  Section~\ref{sect:Discussion} provides a discussion of both the approach and the results.  Finally, Section~\ref{sect:Conclusions} provides conclusions on this research.

\section{General Multiple-Phase Optimal Control Problem \label{sect:multi-phase}}

Without loss of generality, consider the following general multiple-phase optimal control problem, where each phase is defined on the interval $\tau\in[-1,+1]$.
First, let $p \in \{1,\ldots,P\}$ be the phase number where $P$ is the total number of phases.
Determine the state $\m{y}^{(p)}(\tau)\in\bb{R}^{1~\times~n_y^{(p)}}$,
the control $\m{u}^{(p)}(\tau)\in\bb{R}^{1~\times~n_u^{(p)}}$,
the integrals $\m{q}^{(p)}\in\bb{R}^{1~\times~n_q^{(p)}}$, 
the start times $t_0^{(p)} \in\bb{R}$,
and the terminus times $t_f^{(p)} \in\bb{R}$ in all phases $p \in \{1,\dots,P\}$,
along with the static parameters $\m{s} \in\bb{R}^{1~\times~n_s}$ that minimize the objective functional
\begin{equation}\label{eq:multi-cost}
  \C{J}=\C{\phi}\left(\m{e}^{(1)}, \ldots, \m{e}^{(P)}, \m{s}\right)~,
\end{equation}
subject to the dynamic constraints
\begin{equation}\label{eq:multi-dyn}
\frac{dy}{d\tau} \equiv \dot{\m{y}}^{(p)} =
\frac{t_f^{(p)}-t_0^{(p)}}{2}\m{a}^{(p)}\left(\m{y}^{(p)}(\tau),\m{u}^{(p)}(\tau), t^{(p)}(\tau,t_0^{(p)},t_f^{(p)}), \m{s}\right)
~,\quad
\forall p \in \{1,\ldots,P\}~,
\end{equation}
the event constraints
\begin{equation}\label{eq:multi-event}
\m{b}_{\min} \leq \m{b}\left(\m{e}^{(1)}, \ldots, \m{e}^{(P)}, \m{s}\right) \leq \m{b}_{\max}~,
\end{equation}
the inequality path constraints
\begin{equation}\label{eq:multi-path}
\m{c}_{\min}^{(p)} \leq \m{c}^{(p)}\left(\m{y}^{(p)}(\tau),\m{u}^{(p)}(\tau), t^{(p)}(\tau,t_0^{(p)},t_f^{(p)}), \m{s}\right)\leq \m{c}_{\max}^{(p)}
~,\quad
\forall p \in \{1,\ldots,P\}~,
\end{equation}
the integral constraints
\begin{equation}\label{eq:multi-integral}
\m{q}_{\min}^{(p)} \leq \m{q}^{(p)} \leq \m{q}_{\max}^{(p)}
~,\quad
\forall p \in \{1,\ldots,P\}~,
\end{equation}
and the static parameter constraints
\begin{equation}\label{eq:multi-static}
\m{s}_{\min}  \leq \m{s} \leq \m{s}_{\max}~,
\end{equation}
where 
\begin{equation}\label{eq:multi-endpoint-vec}
\m{e}^{(p)} = \left[ \m{y}^{(p)}(-1), t_0^{(p)}, \m{y}^{(p)}(+1), t_f^{(p)}, \m{q}^{(p)} \right]
~,\quad
\forall p \in \{1,\ldots,P\}~,
\end{equation}
\begin{equation}
t^{(p)} \equiv t^{(p)}(\tau,t_0^{(p)},t_f^{(p)}) = \frac{t_f^{(p)}-t_0^{(p)}}{2}\tau + \frac{t_f^{(p)}+t_0^{(p)}}{2}~,
\end{equation}
and the vector components of the integral $\m{q}^{(p)}$ are defined as
\begin{equation}\label{eq:multi-approx}
\begin{array}{rcl}
q_j^{(p)} & = & \displaystyle \frac{t_f^{(p)}-t_0^{(p)}}{2}\int_{-1}^{+1} g_j^{(p)} \left(\m{y}^{(p)}(\tau),\m{u}^{(p)}(\tau), t^{(p)}(\tau,t_0^{(p)},t_f^{(p)}), \m{s}\right) dt~,
\\
 & & \forall j \in \{1,\ldots,n_q^{(p)}\}~,~
\forall p \in \{1,\ldots,P\}~.
\end{array}
\end{equation}
It is noted that the event constraints of Eq.~\eqref{eq:multi-event} contain functions which can relate information at the start and/or terminus of a phase.

\section{Legendre-Gauss-Radau Collocation \label{sect:LGR}}

As stated at the outset, the objective of this research is to compare the efficiency of various derivative approximation methods for use in solving a nonlinear programming problem (NLP) arising from the direct collocation of an optimal control problem.  In order to make the analysis tractable, a particular direct collocation method must be chosen.  While in principle any collocation method can be used to approximate the optimal control problem given in Section \ref{sect:multi-phase}, in this research the basis of comparison will be a recently developed adaptive Gaussian quadrature method called the Legendre-Gauss-Radau (LGR) collocation method \cite{Garg1,Garg2,Garg3,Kameswaran1,Patterson2015,Liu2015,Liu2018}.  It is noted that the NLP arising from the LGR collocation method has an elegant structure.  In addition, the LGR collocation method has a well established convergence theory as described in Refs.~\cite{HagerHouRao15a,HagerHouRao16a,HagerLiuMohapatraWangRao18a}.

In the context of this research, a multiple-interval form of the LGR collocation method is chosen.  In the multiple-interval LGR collocation method, for each phase $p$ of the optimal control problem (where the phase number $p\in\{1,\ldots,P\}$ has been omitted in order to improve clarity of the description of the method), the time interval $\tau\in[-1,+1]$ is divided into $K$ mesh intervals, $\C{S}_k=[T_{k-1},T_k]\subseteq [-1,+1],\;k\in\{1,\ldots,K\}$ such that
\begin{equation}
\bigcup_{k=1}^K\C{S}_k=[-1,+1]
~,
\quad\bigcap_{k=1}^K\C{S}_k=\{T_1,\ldots,T_{K-1}\}~,
\end{equation}
and $-1=T_0<T_1<\ldots<T_{K-1}<T_K=+1$.  For each mesh interval, the LGR points used for collocation are defined in the domain of $[T_{k-1},T_k]$ for $k\in\{1,\ldots,K\}$.  The state of the continuous optimal control problem is then approximated in mesh interval $\C{S}_k,\;k\in\{1,\ldots,K\}$, as 
\begin{equation}\label{eq:LGR-state-approximation}
\m{y}^{(k)}(\tau)  \approx \m{Y}^{(k)}(\tau) = \sum_{j=1}^{N_k+1} \m{Y}_{j}^{(k)}
\ell_{j}^{(k)}(\tau)~, \quad  \ell_{j}^{(k)}(\tau) = \prod_{\stackrel{l=1}{l\neq j}}^{N_k+1}\frac{\tau-\tau_{l}^{(k)}}{\tau_{j}^{(k)}-\tau_{l}^{(k)}}~, 
\end{equation}  
where
$\ell_{j}^{(k)}(\tau)$ for $ j\in\{1,\ldots,N_k+1\}$ is a basis of Lagrange polynomials on $\C{S}_k$, $\left(\tau_1^{(k)},\ldots,\tau_{N_k}^{(k)}\right)$ are the set of $N_k$ Legendre-Gauss-Radau (LGR) \cite{Abramowitz1} collocation points in the interval $[T_{k-1},T_k)$, $\tau_{N_k+1}^{(k)}=T_k$ is a non-collocated support point, and $\m{Y}_{j}^{(k)} \equiv  \m{Y}^{(k)}(\tau_j^{(k)})$.
Differentiating
$\m{Y}^{(k)}(\tau)$ in Eq.~(\ref{eq:LGR-state-approximation}) with
respect to $\tau$ gives
\begin{equation}\label{eq:LGR-diff-state-approximation}
  \frac{d\m{Y}^{(k)}(\tau)}{d\tau} = \sum_{j=1}^{N_k+1}\m{Y}_{j}^{(k)}\frac{d\ell_j^{(k)}(\tau)}{d\tau}~.
\end{equation}
The dynamics are then approximated at the $N_k$ LGR points in mesh
interval $k\in\{1,\ldots,K\}$ as
\begin{equation}\label{eq:LGR-defect}
 \sum_{j=1}^{N_k+1}D_{ij}^{(k)} \m{Y}_j^{(k)} = \frac{t_f-t_0}{2}\m{a}\left(\m{Y}_i^{(k)},\m{U}_i^{(k)},t (\tau_i^{(k)},t_0,t_f),\m{s}\right)
~,\quad
\forall i \in \{1,\ldots,N_k\}~,
\end{equation}
where 
\begin{equation*}
  D_{ij}^{(k)} = \frac{d\ell_j^{(k)}(\tau_i^{(k)})}{d\tau}~,\quad
   \forall i \in \{1,\ldots,N_k\}~,~\forall j \in \{1,\ldots,N_k+1\}~,
\end{equation*}
are the elements of the $N_k\times (N_k+1)$ {\em Legendre-Gauss-Radau differentiation matrix} \cite{Garg1} in mesh interval $\C{S}_k$, $\;k\in\{1,\ldots,K\}$, and $\m{U}_i^{(k)}$ is the approximated control at the $i^{th}$ collocation point in mesh interval $\C{S}_k$.
Finally, introducing the previously omitted phase number, the $P$ phases of the problem are linked together by the event constraints
\begin{equation}\label{eq:LGR-event}
  \m{b}_{\min} \leq \m{b}\left(\m{E}^{(1)}, \ldots, \m{E}^{(P)}, \m{s}\right) \leq \m{b}_{\max}~,
\end{equation}
where $\m{E}^{(p)}$ is the endpoint approximation vector in phase $p$ defined as
\begin{equation}\label{eq:NLP-endpoint-vec}
\m{E}^{(p)} = \left[ \m{Y}_{1}^{(p)}, t_0^{(p)}, \m{Y}_{N^{(p)}+1}^{(p)}, t_f^{(p)}, \m{Q}^{(p)} \right]~,
\end{equation}
such that $N^{(p)}$ is the total number of collocation points used in phase $p$ given by, 
\begin{equation}
N^{(p)} = \sum_{k=1}^{K^{(p)}} N_k^{(p)}~,
\end{equation}
and $\m{Q}^{(p)}\in\mathbb{R}^{1~\times~n_q^{(p)}}$ is the integral approximation vector in phase $p$.

The aforementioned LGR approximation of the continuous optimal control problem leads to the following NLP for a $P$-phase optimal control.  Minimize the objective function
\begin{equation}\label{eq:NLP-cost}
  \C{J}=\C{\phi}\left(\m{E}^{(1)}, \ldots, \m{E}^{(P)}, \m{s}\right)~,
\end{equation}
subject to the defect constraints
\begin{equation}\label{eq:NLP-defect}
\m{\Delta}^{(p)} =
\m{D}^{(p)}\m{Y}^{(p)} - \frac{t_f^{(p)}-t_0^{(p)}}{2}\m{A}^{(p)}
=\m{0}
~,\quad
\forall p \in \{1,\ldots,P\}~,
\end{equation}
the path constraints
\begin{equation}\label{eq:NLP-path}
\m{c}_{\min}^{(p)} \leq
\m{C}_{i}^{(p)}
\leq \m{c}_{\max}^{(p)}
~,\quad
\forall  i \in \{1,\ldots,N^{(p)}\}~,~\forall p \in \{1,\ldots,P\}~,
\end{equation}
the event constraints
\begin{equation}\label{eq:NLP-event}
  \m{b}_{\min} \leq \m{b}\left(\m{E}^{(1)}, \ldots, \m{E}^{(P)}, \m{s}\right) \leq \m{b}_{\max}~,
\end{equation}
the integral constraints
\begin{equation}\label{eq:NLP-integral}
\m{q}_{\min}^{(p)} \leq \m{Q}^{(p)} \leq \m{q}_{\max}^{(p)}
~,\quad
\forall p \in \{1,\ldots,P\}~,
\end{equation}
the static parameter constraints
\begin{equation}\label{eq:NLP-static}
\m{s}_{\min}  \leq \m{s} \leq \m{s}_{\max}~,
\end{equation}
and integral approximation constraints
\begin{equation}\label{eq:NLP-approx}
\g{\rho}^{(p)} =
\m{Q}^{(p)} - \frac{t_f^{(p)}-t_0^{(p)}}{2}\begin{bmatrix}
\m{w}^{(p)}
\end{bmatrix}^{\tr}
\m{G}^{(p)}
= \m{0}
~,\quad
\forall p\in\{1,\ldots,P\}~,
\end{equation}
where
\begin{equation}\label{eq:NLP-A-def}
\m{A}^{(p)} = \begin{bmatrix}
\m{a}^{(p)}\left(\m{Y}_{1}^{(p)},\m{U}_{1}^{(p)},t_1^{(p)},\m{s}\right)\\
\vdots\\
\m{a}^{(p)}\left(\m{Y}_{N^{(p)}}^{(p)},\m{U}_{N^{(p)}}^{(p)},t_{N^{(p)}}^{(p)},\m{s}\right)
\end{bmatrix}
\in \mathbb{R}^{N^{(p)}~\times~n_y^{(p)}}~,
\end{equation}
\begin{equation}\label{eq:NLP-C-def}
\m{C}^{(p)} = \begin{bmatrix}
\m{c}^{(p)}\left(\m{Y}_{1}^{(p)},\m{U}_{1}^{(p)},t_1^{(p)},\m{s}\right)\\
\vdots\\
\m{c}^{(p)}\left(\m{Y}_{N^{(p)}}^{(p)},\m{U}_{N^{(p)}}^{(p)},t_{N^{(p)}}^{(p)},\m{s}\right)
\end{bmatrix}
\in \mathbb{R}^{N^{(p)}~\times~n_c^{(p)}}~,
\end{equation}
\begin{equation}\label{eq:NLP-G-def}
\m{G}^{(p)} = \begin{bmatrix}
\m{g}^{(p)}\left(\m{Y}_{1}^{(p)},\m{U}_{1}^{(p)},t_1^{(p)},\m{s}\right)\\
\vdots\\
\m{g}^{(p)}\left(\m{Y}_{N^{(p)}}^{(p)},\m{U}_{N^{(p)}}^{(p)},t_{N^{(p)}}^{(p)},\m{s}\right)
\end{bmatrix}
\in \mathbb{R}^{N^{(p)}~\times~n_q^{(p)}}~,
\end{equation}
$\m{D}^{(p)} \in \mathbb{R}^{N^{(p)}~\times~[N^{(p)}+1]}$ is the LGR differentiation matrix in phase $p \in \{1,\ldots,P\}$, and $\m{w}^{(p)} \in \mathbb{R}^{N^{(p)}~\times~1}$ are the LGR weights at each node in phase $p$.  It is noted that $\m{a}^{(p)} \in \mathbb{R}^{1~\times~n_y^{(p)}}$, $\m{c}^{(p)} \in \mathbb{R}^{1~\times~n_c^{(p)}}$, and $\m{g}^{(p)} \in \mathbb{R}^{1~\times~n_q^{(p)}}$ correspond, respectively, to the functions that define the right-hand side of the dynamics, the path constraints, and the integrands in phase $p\in\{1,\ldots,P\}$, where $n_y^{(p)}$, $n_c^{(p)}$, and $n_q^{(p)}$ are, respectively, the number of state components, path constraints, and integral components in phase $p$.  Finally, the state matrix, $\m{Y}^{(p)}$, and the control matrix, $\m{U}^{(p)}$, in phase $p\in\{1,\ldots,P\}$ are formed as 
\begin{equation}\label{eq:NLP-Var-Mat}
\m{Y}^{(p)} = \begin{bmatrix}
\m{Y}^{(p)}_{1}\\
\vdots \\
\m{Y}^{(p)}_{N^{(p)}+1}
\end{bmatrix}
\text{ and }
\m{U}^{(p)} = 
\begin{bmatrix}
\m{U}^{(p)}_{1}  \\
\vdots	\\
\m{U}^{(p)}_{N^{(p)}}
\end{bmatrix}~,
\end{equation}
respectively, where $n_u^{(p)}$ is the number control components in phase $p$.

\section{Nonlinear Programming Problem Arising from LGR Collocation\label{sect:NLP}}

\begin{figure}[ht]
\centering
\hspace{-0.2in}
\subfloat[NLP Constraints Jacobian]{\includegraphics[height=3.45in]{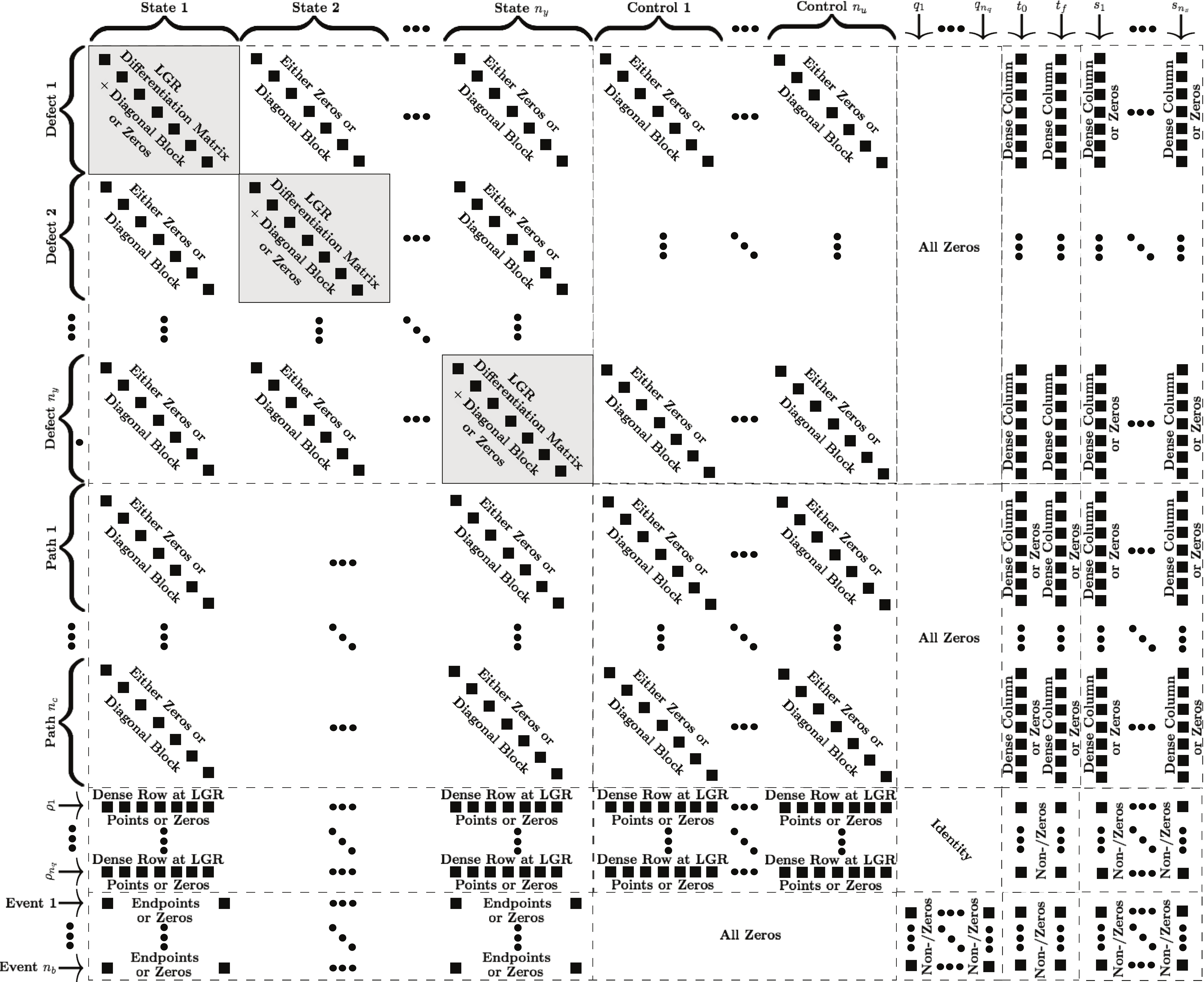}\label{fig:ExampleNLPConJac}}

\subfloat[NLP Lagrangian Hessian]{\includegraphics[height=4.05in]{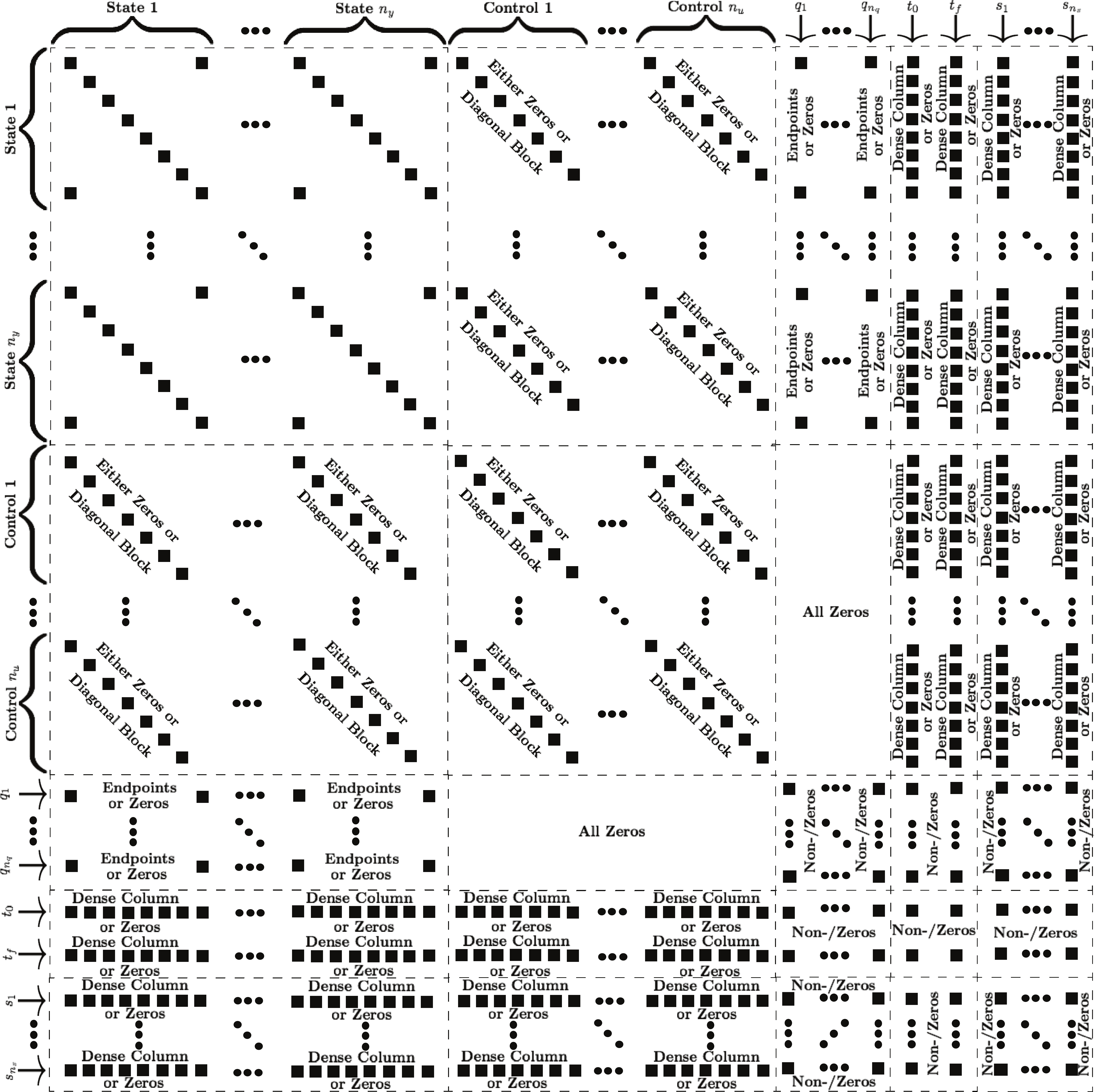}\label{fig:ExampleNLPLagHess}}

\renewcommand{\baselinestretch}{1}\normalsize\normalfont
\caption{Example NLP constraints Jacobian and Lagrangian Hessian sparsity patterns for single-phase optimal control problem with $n_y$ state components, $n_u$ control components, $n_q$ integral components, $n_c$ path constraints, $n_s$ static parameters, and $n_b$ event constraints. \label{fig:ExampleNLPSparsePat}}

\end{figure}

In this section, a brief overview of the resulting nonlinear programming problem (NLP) that arises when discretizing a continuous optimal control problem using Legendre-Gauss-Radau (LGR) collocation is given.  The decision vector of the NLP is comprised of the values of the state at the LGR points plus the final point in each phase, the values of the control at the LGR points in each phase, the initial time of each phase, the final time of each phase, the components of the integral vector in each phase, and any static parameters (which are independent of phase).  The NLP constraints vector is comprised of the defect constraints applied at the LGR points in each phase, the path constraints at the LGR points in each phase, the integral constraints in each phase, and the event constraints (where, in general, the event constraints depend upon information at the start and/or terminus of each phase).  It is noted for completeness that sparsity exploitation as derived in Refs.~\cite{Patterson2012,Patterson2014,Agamawi2018} requires that the first- and second-order partial derivatives of the continuous-time optimal control problem functions be computed at the collocation points \cite{Patterson2012,Patterson2014,Agamawi2018}.  These derivatives are then inserted into the appropriate locations in the NLP constraints Jacobian and Lagrangian Hessian.  Schematics of the LGR collocation NLP derivative matrices are shown in Fig.~\ref{fig:ExampleNLPSparsePat} for a single phase optimal control problem.  It is noted that for the NLP constraints Jacobian, all of the off-diagonal phase blocks relating constraints in phase $i$ to variables in phase $j$ for $i\neq j$ are all zeros.  Similarly, for the NLP Lagrangian Hessian, all of the off-diagonal phase blocks relating variables in phase $i$ to variables in phase $j$ 

\clearpage

\clearpage

\noindent for $i\neq j$ are all zeros except for the variables making up the endpoint vectors which may be related via the objective function or event constraints.  The sparsity patterns shown in Fig.~\ref{fig:ExampleNLPSparsePat} are determined explicitly by identifying the derivative dependencies of the NLP objective and constraints functions with respect to the NLP decision vector variables.  Techniques for identifying the derivative dependencies of the continuous-time optimal control problem functions in terms of the problem variables can significantly increase the sparsity of the NLP derivative matrices by enabling the removal of any elements which are simply zero.  The derivation of the structure of the NLP derivative matrices and the associated sparsity exploitation has been described in detail in Refs.~\cite{Patterson2012,Patterson2014,Agamawi2018} and is beyond the scope of this paper.  

\section{Derivative Approximation Methods \label{sect:deriv-approx}}

As stated at the outset, the objective of this research is to provide a comparison of the use of various Taylor series-based derivative approximation methods for use in direct collocation methods for optimal control.  To this end, this section provides a description of the following three Taylor series-based approximation methods: (1) finite-difference, (2) bicomplex-step, and (3) hyper-dual.  In addition, a description is provided of automatic differentiation.  The computational efficiency of all three Taylor series-based derivative approximation methods are compared to one another and to automatic differentiation.


\subsection{Finite-Difference Methods \label{sect:finite-difference}}

Consider a real function, $f(x):\mathbb{R}\rightarrow\mathbb{R}$ such that $f(x)$ is analytic, infinitely differentiable, and is defined on a domain $\mathcal{D}\subseteq\mathbb{R}$.  By truncating a Taylor series expansion of $f(x_0)$ to the first-order derivative, $f'(x_0)$, the \textit{central finite-difference first-order derivative approximation} can be derived as
\begin{equation}\label{eq:Central-Finite-Difference-First-Order}
f'(x_0) \approx \frac{f(x_0+h)-f(x_0-h)}{2h}~,
\end{equation}
where $h=|x-x_0|$ is the step size taken for the approximation and the truncation error for Eq.~\eqref{eq:Central-Finite-Difference-First-Order} is on the order $\mathcal{O}(h^2)$.  Similarly, the \textit{central finite-difference second-order derivative approximation} can be written as 
\begin{equation}\label{eq:Central-Finite-Difference-Second-Order}
f''(x_0) \approx \frac{f(x_0+h)-2f(x_0)+f(x_0-h)}{h^2}~,
\end{equation}
where the truncation error for Eq.~\eqref{eq:Central-Finite-Difference-Second-Order} is on the order $\mathcal{O}(h^2)$.  Extending the concepts of finite-difference methods to multivariable functions, for a real function $g(x,y):\mathbb{R}\times\mathbb{R}\rightarrow\mathbb{R}$, the second-order mixed partial derivative of $g$ at $x_0,y_0 \in \mathbb{R}$ can be approximated with the central finite-difference method as
\begin{equation}\label{eq:Central-Finite-Difference-Second-Order-Mixed}
\begin{array}{rcl}
\displaystyle \frac{\partial^2 g(x_0,y_0)}{\partial x \partial y} & \approx & \displaystyle \left(g(x_0+h_x,y_0+h_y)-g(x_0+h_x,y_0-h_y)-g(x_0-h_x,y_0+h_y)\right. \\
 & + & \left.g(x_0-h_x,y_0-h_y)\right)/\left(4h_xh_y\right)~,
\end{array}
\end{equation}
where $h_x$ and $h_y$ are the step sizes taken for the independent variables $x_0$ and $y_0$, respectively. The truncation error for Eq.~\eqref{eq:Central-Finite-Difference-Second-Order-Mixed} is on the order $\mathcal{O}(h_x h_y)$.

In addition to truncation error, the derivative approximations in Eqs.~\eqref{eq:Central-Finite-Difference-First-Order} to \eqref{eq:Central-Finite-Difference-Second-Order-Mixed} are all subject to roundoff error due to the differences being taken to compute the approximations.  Due to finite-difference methods being subject to both truncation and roundoff error, when implementing finite-differencing, an appropriate step size, $h$, must be chosen in order to reduce the truncation error while refraining from exacerbating roundoff error.  For the purposes of using finite-differencing as a method for an optimization software, the appropriate step size, $h$, may be computed as,
\begin{equation}\label{eq:Finite-Difference-Step-Size}
h = h_{\mathcal{O}_1}(1+|x_0|)~,
\end{equation}
where $|x_0|$ is the magnitude of the independent variable at the point of interest and the base perturbation size, $h_{\mathcal{O}_1}$, is chosen to be the optimal step size for a function whose input and output are $\approx \mathcal{O}(1)$ as described in
Ref.~\cite{Gill3}.

\subsection{Bicomplex-step Derivative Approximation \label{sect:Bicomplex}}

In this section, the basis of bicomplex-step derivative approximations are derived.  As described in Ref.~\cite{Lantoine1}, an efficient and accurate way to compute derivatives required by the NLP solver is to use the bicomplex-step derivative approximation.  Consider a complex number
\begin{equation}\label{eq:Complex-Number}
c \in \mathbb{C}^1 = x + i_1y~,
\end{equation}
where $x,y \in \mathbb{R} \equiv \mathbb{C}^0$ and $i_1$ designates the imaginary component of the complex number in the complex plane $\mathbb{C}^1$ and has the property $i_1^2 = -1$.  Extending the concept of complex numbers to second-order imaginary numbers, a bicomplex number is defined as
\begin{equation}\label{eq:Bicomplex-Number}
z \in \mathbb{C}^2 = c_{1} + i_2c_{2}~,
\end{equation}
where $c_{1},c_{2} \in \mathbb{C}^1$ and $i_2$ designates an imaginary component distinct from that of the $i_1$ imaginary direction.  The additional imaginary direction, $i_2$, is similar to $i_1$ in that $i_2^2 = -1$.  For bicomplex numbers, $i_1 i_2$ remains $i_1i_2$ (that is, $i_1i_2$ is considered a \textit{bi-imaginary} component distinct from either $i_1$ or $i_2$) and the multiplicative relationship between imaginary number components is commutative (that is $i_1 i_2 = i_2 i_1$).  Using the definitions of complex numbers from Eq.~\eqref{eq:Complex-Number}, Eq.~\eqref{eq:Bicomplex-Number} may be expanded and written in terms of its real number components as
\begin{equation}\label{eq:Bicomplex-Number-Expanded}
z = x_1 + i_1y_1 + i_2x_2 + i_1i_2y_2~.
\end{equation}

Using this definition of bicomplex numbers, consider a bicomplex function $f(z):\mathbb{C}^2\rightarrow\mathbb{C}^2$ that is analytic in the domain $\mathcal{D}\in\mathbb{C}^2$.  Because the functions of interest are real-valued, the values of $z$ can be restricted to the domain $\mathcal{D}\in\mathbb{R}$ such that $f(x):\mathbb{R}\rightarrow\mathbb{R}$ is isomorphic with respect to $f(z),~\forall z\in\mathcal{D}$ (that is, bicomplex numbers with zero imaginary components).  For $z_0\in\mathcal{D}$, a Taylor series expansion may be taken about the real component $x$ of the bicomplex number $z_0$ by taking purely imaginary steps $hi_1$ and $hi_2$ in both the imaginary component directions such that
\begin{equation}\label{eq:Bicomplex-Step}
\begin{array}{rcl}
z_0 & = & x~, \\
z & = & x + i_1h + i_2h~.
\end{array}
\end{equation}
Extending the concept of a Taylor series expansion of a bicomplex function to a multivariable bicomplex function, the following partial derivative approximations can be made for a given analytic function\\
\noindent $f(x,y):\mathbb{R}\times\mathbb{R}\rightarrow\mathbb{R}$:
\begin{equation}\label{eq:Bicomplex-step-Derivative-Approximations}
\begin{array}{l}
\begin{array}{rcl}
\displaystyle \frac{\partial^2f(x,y)}{\partial x^2} &  \approx & \displaystyle \frac{\text{Im}_{1,2}\left[f(x+hi_1+hi_2,y)\right]}{h^2}~,  \vspace{0.05in}\\
\displaystyle \frac{\partial^2f(x,y)}{\partial y^2} &  \approx & \displaystyle \frac{\text{Im}_{1,2}\left[f(x,y+hi_1+hi_2)\right]}{h^2}~,
\vspace{0.05in}\\
\displaystyle \frac{\partial^2f(x,y)}{\partial x \partial y} &  \approx & \displaystyle \frac{\text{Im}_{1,2}\left[f(x+hi_1,y+hi_2)\right]}{h^2}~,  \vspace{0.05in}\\
\displaystyle \frac{\partial f(x,y)}{\partial x} &  \approx & \displaystyle \frac{\text{Im}_{1}\left[f(x+hi_1+hi_2,y)\right]}{h} = \frac{\text{Im}_{2}\left[f(x+hi_1+hi_2,y)\right]}{h}~, \vspace{0.05in}\\
\displaystyle \frac{\partial f(x,y)}{\partial y} &  \approx & \displaystyle \frac{\text{Im}_{1}\left[f(x,y+hi_1+hi_2)\right]}{h} = \frac{\text{Im}_{2}\left[f(x,y+hi_1+hi_2)\right]}{h}~,
\vspace{0.05in}\\
\displaystyle \frac{\partial f(x,y)}{\partial x} &  \approx & \displaystyle \frac{\text{Im}_{1}\left[f(x+hi_1,y+hi_2)\right]}{h}~,  \vspace{0.05in}\\
\displaystyle \frac{\partial f(x,y)}{\partial y} &  \approx & \displaystyle \frac{\text{Im}_{2}\left[f(x+hi_1,y+hi_2)\right]}{h}~,
\end{array}
\end{array}
\end{equation}
where $h>0$ is the step size taken in a given imaginary direction and the Taylor series expansion is truncated to the second-order derivative such that the truncation error is $\C{O}(h^2)$. Thus for any given function $f(x,y)$ evaluated at inputs $x,y \in \mathbb{R}$, the first- and second-order partial derivatives may be obtained by evaluating the function with bicomplex inputs that have the same real parts and additionally the appropriate imaginary direction components so to obtain the desired partial derivatives.  Such an approximation for the first- and second-order derivatives is extremely useful, as it avoids the need to explicitly (or implicitly) derive the differential equations.  Additionally, bicomplex-step derivative approximations do not require taking any differences, thus avoiding roundoff error associated with finite-differencing and allowing the step sizes taken in the imaginary directions to be made arbitrarily small in order to reduce the truncation error.  It is noted that although bicomplex-step derivative approximations avoid roundoff error in the method itself, the bicomplex arithmetic necessary for evaluating the functions are still subject to roundoff error in implementation, thus limiting the step size from being too small \cite{Lantoine1}.  Moreover, the independent nature of the bicomplex-step derivative approximations for first-order derivatives enables  two independent first-order partial derivative approximations to be computed for the same function using a single bicomplex evaluation.  Furthermore, the ease of using the bicomplex-step derivative approximation to simply evaluate the function at the point of interest using bicomplex inputs with appropriate imaginary components can be taken advantage of by using operator overloading in C++.  By defining a bicomplex class and using operator overloading for all elementary functions, a function written in C++ can be easily evaluated with a bicomplex input, and the resulting bicomplex output can be used to determine the desired partial derivatives.
The bicomplex-step derivative approximation thus provides a fast and accurate means of determining the first- and second-order partial derivatives of a given function.

\subsection{Hyper-Dual Derivative Approximation \label{sect:HyperDual}}

In this section, the basis of hyper-dual derivative approximations are derived.  As described in Ref.~\cite{Fike2011}, an efficient and exact way to compute derivatives required by the NLP solver is to use the hyper-dual derivative approximation.  Consider a dual number,
\begin{equation}\label{eq:dual}
d\in\mathbb{D}^1 = x + \epsilon_1 y~,
\end{equation}
where $x,y\in\mathbb{R}\equiv \mathbb{D}^0$ and $\epsilon_1$ designates the imaginary component of the dual number in the dual plane $\mathbb{D}^1$ and has the property $\epsilon_1^2 = 0$ (that is, $\epsilon_1$ is nilpotent).  Extending the concept of a dual number to second-order imaginary numbers, a hyper-dual number is defined as
\begin{equation}\label{eq:hyperdual}
w\in\mathbb{D}^2 = d_1 + \epsilon_2 d_2~,
\end{equation}
where $d_1,d_2\in\mathbb{D}^1$ and $\epsilon_2$ designates an imaginary component distinct from that of the $\epsilon_1$ imaginary direction. The additional imaginary direction, $\epsilon_2$, is similar to $\epsilon_1$ in that $\epsilon_2^2=0$.  The multiplicative relationship between the imaginary directions of the hyper-dual plane is also analogous to how the bicomplex imaginary directions relate as described in Section~\ref{sect:Bicomplex} (that is $\epsilon_1\epsilon_2$ remains $\epsilon_1\epsilon_2$ and the multiplication of components is commutative). 

Using this definition of hyper-dual numbers, consider a hyper-dual function $f(w):\mathbb{D}^2\rightarrow\mathbb{D}^2$ that is analytic in the domain $\mathcal{E}\in\mathbb{D}^2$.  Because the functions of interest are real-valued, values of $w$ can be restricted to the domain $\mathcal{E}\in\mathbb{R}$ such that $f(w):\mathbb{R}\rightarrow\mathbb{R}$ is isomorphic with respect to $f(w),~\forall w\in\mathcal{E}$ (that is, hyper-dual numbers with zero imaginary components).  For $w_0\in\mathcal{E}$, a Taylor series expansion may be taken about the real component $x$ of the hyper-dual number $w_0$ by taking purely imaginary steps $h\epsilon_1$ and $h\epsilon_2$ in both the imaginary component directions such that
\begin{equation}\label{eq:HyperDualInput}
\begin{array}{rcl}
w_0 & = & x~, \\
w & = & x + \epsilon_1h + \epsilon_2h~.
\end{array}
\end{equation}
Extending the concept of a Taylor series expansion of a hyper-dual function to a multivariable hyper-dual function, the following partial derivative approximations can be made for a given analytic function $f(x,y):\mathbb{R}\times\mathbb{R}\rightarrow\mathbb{R}$:
\begin{equation}\label{eq:hyperdual-mixed-approx}
\begin{array}{l}
\begin{array}{rcl}
\displaystyle \frac{\partial^2f(x,y)}{\partial x^2} & = & \displaystyle \frac{\text{Ep}_{1,2}\left[f(x+h\epsilon_1+h\epsilon_2,y)\right]}{h^2}~,  \vspace{0.05in}\\
\displaystyle \frac{\partial^2f(x,y)}{\partial y^2} & = & \displaystyle \frac{\text{Ep}_{1,2}\left[f(x,y+h\epsilon_1+h\epsilon_2)\right]}{h^2}~,
\vspace{0.05in}\\
\displaystyle \frac{\partial^2f(x,y)}{\partial x \partial y} & = & \displaystyle \frac{\text{Ep}_{1,2}\left[f(x+h\epsilon_1,y+h\epsilon_2)\right]}{h^2}~,  \vspace{0.05in}\\
\displaystyle \frac{\partial f(x,y)}{\partial x} & = & \displaystyle \frac{\text{Ep}_{1}\left[f(x+h\epsilon_1+h\epsilon_2,y)\right]}{h} = \frac{\text{Ep}_{2}\left[f(x+h\epsilon_1+h\epsilon_2,y)\right]}{h}~, \vspace{0.05in}\\
\displaystyle \frac{\partial f(x,y)}{\partial y} & = & \displaystyle \frac{\text{Ep}_{1}\left[f(x,y+h\epsilon_1+h\epsilon_2)\right]}{h} = \frac{\text{Ep}_{2}\left[f(x,y+h\epsilon_1+h\epsilon_2)\right]}{h}~,
\vspace{0.05in}\\
\displaystyle \frac{\partial f(x,y)}{\partial x} & = & \displaystyle \frac{\text{Ep}_{1}\left[f(x+h\epsilon_1,y+h\epsilon_2)\right]}{h}~,  \vspace{0.05in}\\
\displaystyle \frac{\partial f(x,y)}{\partial y} & = & \displaystyle \frac{\text{Ep}_{2}\left[f(x+h\epsilon_1,y+h\epsilon_2)\right]}{h}~,
\end{array}
\end{array}
\end{equation}
where $h>0$ is the step size taken in a given imaginary direction and there is no truncation error in any of the approximations due to the nilpotent property of the hyper-dual number imaginary directions.
The partial derivative approximations given in Eq.~\eqref{eq:hyperdual-mixed-approx} are known as \textit{hyper-dual derivative approximations}.  Thus for any given function $f(x,y)$ evaluated at inputs $x,y \in \mathbb{R}$, the first- and second-order partial derivatives may be obtained by evaluating the function with hyper-dual inputs that have the same real parts and additionally the appropriate imaginary direction components so to obtain the desired partial derivatives.  Such an approximation for the first- and second-order derivatives is extremely useful, as it avoids the need to explicitly (or implicitly) derive the differential equations.  Additionally, hyper-dual derivative approximations do not require taking any differences and avoid roundoff error associated with finite-differencing.  Due to the properties of the hyper-dual arithmetic \cite{Fike2011}, the hyper-dual derivative approximations also avoid any susceptibility to roundoff error during function evaluation so as to maintain machine precision for all step sizes.  Moreover, there is no truncation error when using the hyper-dual derivative approximations so the step size may be chosen arbitrarily.  Similar to the bicomplex-step derivative approximation, the independent nature of the hyper-dual derivative approximations for first-order derivatives enables two independent first-order partial derivative approximations to be computed for the same function using a single hyper-dual evaluation, and C++ operator overloading may be employed to allow easy evaluation of functions with hyper-dual inputs.  The hyper-dual derivative approximation thus provides a fast and accurate means of determining the first- and second-order partial derivatives of a given function.

\subsection{Automatic Differentiation \label{sect:AlgorithmicDiff}}

In this section, the basis of automatic differentiation is discussed.  As described in Ref.~\cite{Martins2}, automatic (algorithmic) differentiation may be derived from the unifying chain rule and supplies numerical evaluations of the derivative for a defined computer program by decomposing the program into a sequence of elementary function operations and applying the calculus chain rule algorithmically through the computer \cite{Griewank0}.  The process of automatic differentiation is described in detail in Ref.~\cite{Griewank0}, and is beyond the scope of this paper.  It is noted, however, that the Taylor series-based derivative methods described in Sections~\ref{sect:finite-difference} through \ref{sect:HyperDual} are compared with the open source automatic differentiation software ADOL-C \cite{Griewank1,Walther1}.

\section{Comparison of Various Derivative Approximation Methods \label{sect:DerivativeComparison}}

The purpose of this section is to highlight the important aspects and differences of implementing the derivative estimation methods described in Section~\ref{sect:deriv-approx}.

\subsection{Derivative Approximation Error \label{subsect:Deriv-Approx-Error}}

The estimates of the truncation error for the central finite-difference, bicomplex-step, and hyper-dual derivative approximation methods derived in Sections~\ref{sect:finite-difference} -- \ref{sect:HyperDual}, respectively, are based on a Taylor series expansion of an analytic function.  Both the central finite-difference and the bicomplex-step derivative approximations have a truncation error that is $\mathcal{O}(h^2)$, where $h$ is the step size used for the derivative approximation.  On the other hand, the hyper-dual derivative approximation has no truncation error.  Moreover, finite-difference methods are inherently subject to roundoff errors due to the difference quotient used in the method.  Although the bicomplex-step derivative approximation does not inherently appear to suffer from roundoff error, the evaluation of the bicomplex arithmetic for complicated functions (for example, logarithmic, power, or inverse trigonometric functions) can lead to roundoff error in the derivative approximation.  On the other hand, the hyper-dual derivative approximation avoids roundoff error altogether, as the method itself requires no differences to be taken, and the values computed for the first- and second-order derivative levels are accurate to machine precision regardless of the step size used for the perturbations in the imaginary directions.

\begin{figure}[ht]
\centering
\subfloat[First-Order Derivative Approximation.\label{fig:ExampleFuncFirstDerErr}]{\includegraphics[width=3.2in]{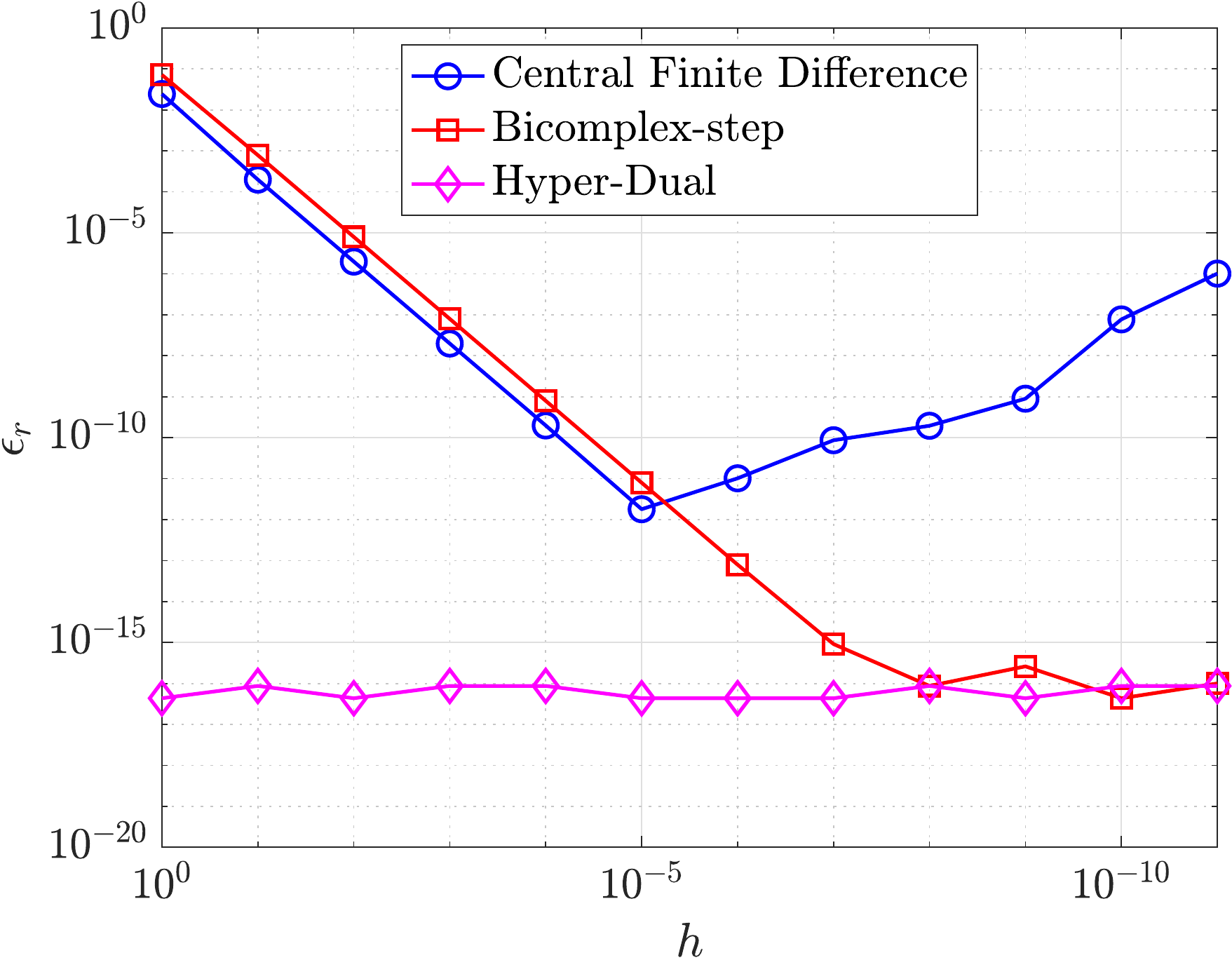}}~\subfloat[Second-Order Derivative Approximation.\label{fig:ExampleFuncSecondDerErr}]{\includegraphics[width=3.2in]{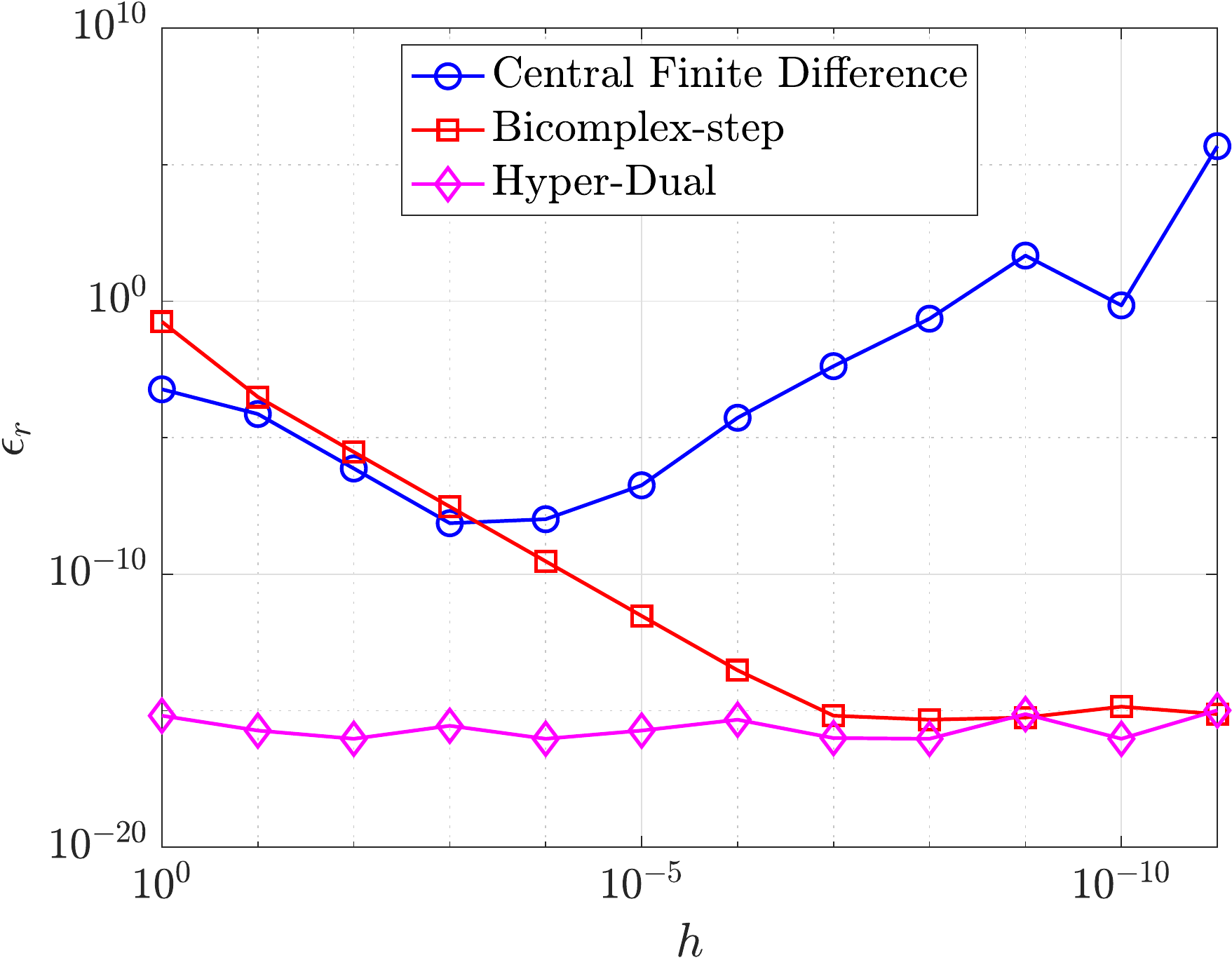}}

\renewcommand{\baselinestretch}{1}\normalsize\normalfont
\caption{Comparison of relative error of derivative approximations obtained using central finite-difference, bicomplex-step, and hyper-dual methods for example function in Eq.~\eqref{eq:Example-Error-Function}. \label{fig:ExampleFuncDerErr}}

\end{figure}

Suppose now that the relative error in a quantity $d$ is defined as
\begin{equation}\label{relative-error}
  \epsilon_r = \frac{|d-\hat{d}|}{1+|d|}~,
\end{equation}
where $\hat{d}$ is the approximation of $d$.  Using the definition of the relative error given in Eq.~\eqref{relative-error}, the accuracy of the three Taylor series-based derivative approximations is shown in Fig.~\ref{fig:ExampleFuncDerErr} for the example function
\begin{equation}\label{eq:Example-Error-Function}
  f(x) = \frac{x^2}{\sin(x)+x\exp(x)}~,
\end{equation}
evaluated at $x=0.5$.  It is noted for completeness that the analytic first- and second-order derivatives of Eq.~\eqref{eq:Example-Error-Function} are given, respectively, as
\begin{equation}\label{eq:Example-Error-Function-First-Derivative}
  \begin{array}{lcl}
    \displaystyle \frac{df(x)}{dx} & = & \displaystyle \frac{x\left(2\sin(x)-x\cos(x)+\left(x-x^2\right)\exp(x)\right)}{\left(\sin(x)+x\exp(x)\right)^2}~, \vspace{3pt} \\
    \displaystyle \frac{d^2f(x)}{dx^2} & = & \left( \left( x^2+2\right)\sin^2(x)+\left(-4x\cos(x)-6x^2\exp(x)\right)\sin(x)+2x^2\cos^2(x) \right. \\
 & + & \left. 4x^3\exp(x)\cos(x)+\left(x^4-2x^3\right)\exp(2x)\right)/\left(\sin(x)+x\exp(x)\right)^3~.
  \end{array}
\end{equation}
It is seen from Fig.~\ref{fig:ExampleFuncDerErr} that the derivative approximation error using the central finite-difference method decreases until $h\approx 10^{-5}$ and then starts to increase for $h<10^{-5}$.  Next, it is seen from Fig.~\ref{fig:ExampleFuncDerErr} that the bicomplex-step derivative approximation error decreases as $h$ decreases until the error reaches near machine precision (with a value $\epsilon_r\approx 10^{-15}$).  Finally, it is seen from Fig.~\ref{fig:ExampleFuncDerErr} that the hyper-dual derivative approximation error remains near machine precision regardless of the value of $h$.  It is noted that the derivative approximation error for automatic differentiation is not shown in Fig.~\ref{fig:ExampleFuncDerErr} because automatic differentiation (see Section~\ref{sect:AlgorithmicDiff}) employs the calculus chain rule and, thus, is accurate to machine precision.

\subsection{Computational Efficiency Expectations \label{subsect:Deriv-Approx-Computational-Efficiency}}

As described in Sections~\ref{sect:finite-difference} -- \ref{sect:HyperDual}, for first- and second-order derivative approximations, central finite-differencing requires two and four real function evaluations, respectively, while both the bicomplex-step and hyper-dual derivative approximations require just one function evaluation for both first- and second-derivatives.  Although the required number of function evaluations may make it appear as if the bicomplex-step and hyper-dual derivative approximations are more computationally efficient, neither of these two methods is guaranteed to be more computationally efficient than central finite-differencing.  In particular, despite requiring more function evaluations, central finite-differencing only requires real arithmetic, while the bicomplex-step and hyper-dual derivative approximations require bicomplex and hyper-dual number arithmetic, respectively.  Specifically, both bicomplex and hyper-dual numbers consist of four real components that lie along distinct directions.  Therefore, a function evaluation using either a bicomplex or hyper-dual number requires at least four times as many operations as the function evaluation of a real number.  For example, a single addition or subtraction of two bicomplex numbers requires four real number operations.  Furthermore, functions such as inverse trigonometric functions and logarithmic functions require an even larger number of operations on real numbers when using either bicomplex or hyper-dual arithmetic.  It is noted that hyper-dual derivative approximations require between one and $3.5$ times more operations than central finite-differencing \cite{Fike2011}, and bicomplex-step derivative approximations may require even more operations due to the complexity of the bicomplex arithmetic \cite{Lantoine1}.  Consequently, although central finite-differencing requires multiple real function evaluations to compute the derivative approximation, it may still be the case that central finite-differencing will require fewer real arithmetic function evaluations when compared with a single hyper-dual or bicomplex function evaluation.  It is noted that the computational efficiency of the bicomplex-step and hyper-dual methods can be improved significantly by taking advantage of the ability to compute simultaneous partial derivatives of the same function using a single bicomplex or hyper-dual function evaluation.  In fact, whenever possible, it is important to use a single bicomplex or hyper-dual function evaluation to simultaneously estimate multiple derivatives so that the number of function evaluations is minimized.  The performance of the three Taylor series-based methods relative to one another in regards to computational efficiency will thus depend on the overall complexity of the function whose partial derivatives are of interest, as well as the number of partial derivatives needed for the function.

Finally, although the automatic differentiation described in Section~\ref{sect:AlgorithmicDiff} utilizes real number arithmetic, the number of operations performed to compute the derivative estimation may be substantial as a result of the tracing process used to tape the record of the function evaluation that is of interest that is required when using automatic differentiation software such as ADOL-C \cite{Griewank1,Walther1}.  Due to the complex nature of optimal control problems, the flow control of the computer programs defining the optimal control problem functions are often conditionally-based on the function inputs (for example, as a result of piecewise models or interpolation of tabular data).  Any such conditional flow control requires the tracing process of the automatic differentiation method to be invoked each time a derivative is needed, as the flow control must be properly recorded for any given input value, thus leading to a substantial amount of computation time to supply derivative estimates \cite{Walther1}.  Consequently, even though the derivative obtained using automatic differentiation is of very high quality, the repeated retracing of the function to compute the derivative significantly decreases computational efficiency.

\subsection{Identification of Derivative Dependencies \label{subsect:Deriv-Approx-Dependencies}}

In order to maximize computational efficiency when solving an NLP arising from a direct collocation method, it is important that the fewest number of derivatives be computed.  In particular, many of the functions of the original continuous optimal control problem (for example, the vector fields that define the right-hand side of the differential equations, the path constraints, and the integrands) are functions of only some of the components of the state and control.  Consequently, the process of exploiting the sparse structure of the NLP leads to many zero blocks in the constraints Jacobian and Lagrangian Hessian because of the independence of these functions with respect to certain variables.  Moreover, eliminating these derivatives from the constraints Jacobian and Lagrangian Hessian can significantly improve the computational efficiency with which the NLP is solved.  The aforementioned discussion makes it important and useful to determine the NLP derivative dependencies based on the dependencies of the optimal control problem functions.

As it turns out, the bicomplex-step and hyper-dual derivative approximation methods described in Sections~\ref{sect:Bicomplex} and \ref{sect:HyperDual} have a key property that the the imaginary components for either of these methods have an independent nature that enables an accurate determination of the optimal control problem derivative dependencies which, in turn, enables computing the fewest derivatives required for use by the NLP solver.  As stated, determining these derivative dependencies leads to full exploitation of the sparsity in the NLP as described in Section~\ref{sect:NLP}.  On the other hand, central finite-differencing as described in Section~\ref{sect:finite-difference} cannot easily be utilized as a tool for identifying derivative dependencies because central finite-differencing is highly susceptible to both truncation and roundoff error.  In particular, because the step size used for central finite-differencing cannot be made arbitrarily small, an approximation error is always present.  Thus, when trying to identify second-order derivative sparsity patterns, the evaluation of mixed partial derivatives may be nonzero even when the partial derivative itself is zero.

Because central finite-differencing cannot be employed effectively for the determination of derivative dependencies, other approaches need to be employed if central finite-differencing is chosen as the derivative estimation method when solving the NLP.  One possibility for determining derivative dependencies when using central finite-differences is to use a technique that employs NaN (not a number) or Inf (infinity) propagation.  Methods that employ NaN (or Inf) propagation involve evaluating a function of interest for variable dependence by setting a variable of interest equal to NaN (or Inf) and seeing if the output of the function of interest returns NaN (or Inf), which would indicate the function of interest is dependent on the variable of interest.  It is noted, however, that NaN and Inf approaches are limited in that they can only identify first-order derivative dependencies.  As a result, NaN and Inf propagation methods lead to an over-estimation of second-order derivative dependencies where some derivatives, which may actually be zero, are estimated to be nonzero.  This over-estimation of the second-order derivative dependencies leads to a denser NLP Lagrangian Hessian which can significantly decrease computational efficiency when solving the NLP.  On the other hand, using either the bicomplex-step or hyper-dual derivative approximation methods, the sparsity patterns of the NLP derivative matrices can be determined exactly which, in turn, can significantly improve computational efficiency when solving the NLP.  It is noted that the automatic differentiation as discussed in Section~\ref{sect:AlgorithmicDiff} may also be utilized to obtain an exact first- and second-order NLP derivative matrices sparsity pattern because the obtained derivative estimates are accurate to machine precision.

\section{Examples \label{sect:Examples}}

In this section, the derivative estimation methods described in Section~\ref{sect:deriv-approx} are compared to one another in terms of their effectiveness as methods for use when solving the nonlinear programming problem (NLP) resulting from the transcription of the continuous optimal control problem using Legendre-Gauss-Radau (LGR) collocation, as developed in Sections~\ref{sect:multi-phase} through \ref{sect:NLP}.  The first example is the free-flying robot optimal control problem taken from Ref.~\cite{Sakawa}.  The second example is the minimum time-to-climb of a supersonic aircraft optimal control problem taken from Ref.~\cite{Bryson1}.  The third example problem is the space station attitude control optimal control problem taken from Ref.~\cite{Betts3}.  All three examples demonstrate that the more accurate derivative estimates of the bicomplex-step, hyper-dual, and automatic differentiation methods as described in Section~\ref{subsect:Deriv-Approx-Error} can significantly reduce the number of NLP solver iterations required to solve the NLP as compared to when using a central finite-difference method.  Note, however, that the results also show the improved computational efficiency obtained using either the bicomplex-step or hyper-dual methods relative to central finite-differencing because the NLP solver requires fewer NLP iterations to converge to a solution using the bicomplex-step and hyper-dual methods.  Furthermore, even though the derivative obtained using automatic differentiation is accurate to machine precision, the lower computational efficiency of the automatically computed derivative results in a significant decrease in computational efficiency when solving the NLP.  Finally, the increase in computational efficiency provided by exactly identifying the first- and second-order derivative dependencies of the continuous-time functions as discussed in Section~\ref{subsect:Deriv-Approx-Dependencies} is highlighted by all three examples.

The following terminology is used in each example. First, $K$ denotes the number of mesh intervals used to discretize the continuous optimal control problem, where the number of collocation points used in each interval $k\in\{1,\ldots,K\}$ is $N_k$.  Furthermore, the notation $OC$, $EC$, $BC$, $HD$ and $AD$ is used to denote the over-estimated sparsity central finite-difference, exact sparsity central finite-difference, bicomplex-step, hyper-dual, and automatic differentiation methods, respectively.  The over-estimated sparsity central finite-difference ($OC$) method refers to using central finite-differencing with the NLP derivative matrices sparsity pattern obtained from identifying the derivative dependencies using NaN propagation, and the exact sparsity central finite-difference ($EC$) method refers to using central finite-differencing with the NLP derivative matrices sparsity pattern obtained by using the hyper-dual derivative approximations to identify the derivative dependencies, as discussed in Section~\ref{subsect:Deriv-Approx-Dependencies}.  It is noted that the bicomplex-step ($BC$), hyper-dual ($HD$), and automatic differentiation ($AD$) methods utilize the same NLP derivative matrices sparsity pattern as the exact sparsity central finite-difference ($EC$) method.  Moreover, $I$, $\C{T}$, and $\Phi$ denote the number of NLP solver iterations to converge, the total computation time required to solve the NLP, and the average computation time expended per iteration to compute the derivative approximations.  Finally, the use of $\Gamma_{\xi}^{\mathcal{A}}$ for $\xi \in \{I,\C{T},\Phi \}$, $\mathcal{A}\in\{EC,BC,HD,AD\}$, denotes the percent reduction between the value of $\xi$ required using the $\mathcal{A}$ method versus the $OC$ method, as given by
\begin{equation}\label{eq:percent-reduction}
\Gamma_{\xi}^{\mathcal{A}} = \frac{\xi_{OC} - \xi_{\mathcal{A}}}{\xi_{OC}}\times 100~,\quad \xi \in \{I,\C{T},\Phi \}~,\quad\mathcal{A}\in\{EC,BC,HD,AD\}~.
\end{equation}
It is noted that the $OC$ method is used as a baseline of performance comparison in order to highlight the reduced (or increased) computational expense acquired when utilizing the other methods.

All results shown in this paper were obtained using the C++ optimal control software $\mathbb{CGPOPS}$ \cite{Agamawi2018} using the NLP solver IPOPT \cite{Biegler2}.  The NLP solver was set to an optimality tolerance of $10^{-7}$ for all three examples.   All first- and second-order derivatives for the NLP solver were obtained using the derivative estimation methods described in Section~\ref{sect:deriv-approx}.  All computations were performed on a 2.9 GHz Intel Core i7 MacBook Pro running MAC OS-X version 10.13.6 (High Sierra) with 16GB 2133MHz LPDDR3 of RAM.  C++ files were compiled using Apple LLVM version 9.1.0 (clang-1000.10.44.2).

\subsection{Example 1:  Free-Flying Robot Problem \label{subsect:freeFlyingRobot}}

Consider the following optimal control problem taken from Ref.~\cite{Betts3} and \cite{Sakawa}.  Minimize the cost functional  
\begin{equation}\label{eq:freeFlyingCost}
  \C{J} = \int_{0}^{t_f}(u_1(t)+u_2(t)+u_3(t)+u_4(t))~dt~,
\end{equation}
subject to the dynamic constraints
\begin{equation}\label{eq:freeFlyingDynamics}
 \begin{array}{lclclcl}
   \dot{x}(t) & = & v_x(t)~, & & \dot{y}(t) & = & v_y(t)~, \\
   \dot{v_x}(t) & = & (F_1(t)+F_2(t))\cos(\theta(t))~, & & \dot{v_y}(t) & = & (F_1(t)+F_2(t))\sin(\theta(t))~, \\
   \dot{\theta}(t) & = & \omega(t)~, & & \dot{\omega}(t) & = & \alpha F_1(t)-\beta F_2(t)~, \\
\end{array}
\end{equation}
the control inequality constraints
\begin{equation}\label{eq:freeFlyingControlConstraint}
\begin{array}{lcl}
 0 \leq u_i(t) \leq 1000~, \quad (i=1,2,3,4)~, & & \quad F_i(t) \leq 1~,~\quad (i=1,2)~,
\end{array}  
\end{equation}
and the boundary conditions
\begin{equation} \label{eq:freeFlyingBCs}
\begin{array}{lclclcl}
 x(0) = -10~, & & x(t_f) = 0~, & &
 y(0) = -10~, & & y(t_f) = 0~, \\
 v_x(0) = 0~, & & v_x(t_f) = 0~, & &
 v_y(0) = 0~, & & v_y(t_f) = 0~, \\
 \theta(0) = \frac{\pi}{2}~, & & \theta(t_f) = 0~, & &
 \omega(0) = 0~, & & \omega(t_f) = 0~,
\end{array}
\end{equation} 
where 
\begin{equation}
\label{eq:freeFlyingAux}
 \begin{array}{lclclcl}
 F_1(t) = u_1(t)-u_2(t)~, & & F_2(t) = u_3(t)-u_4(t)~,& &
 \alpha = 0.2~, & & \beta = 0.2~.
 \end{array}
\end{equation}
The optimal control problem described by Eqs.~\eqref{eq:freeFlyingCost} -- \eqref{eq:freeFlyingAux} was approximated as an NLP using LGR collocation (see Section~\ref{sect:LGR}) using $K=(2,4,8,16,32)$ mesh intervals with $N_k=5$ in each mesh interval.  The resulting set of NLPs for the different values of $K$ were then solved using an over-estimated and exact sparsity with central finite-difference,  and using exact sparsity patterns with bicomplex-step, hyper-dual, and automatic differentiation.  The solutions obtained in all cases are essentially identical and converge to the exact solution given in Ref.~\cite{Betts3} as $K$ increases (see pages $328$ to $329$ of Ref.~\cite{Betts3}).  Comparisons of the computational efficiency in terms of $I$, $\C{T}$, and $\Phi$ for varying values of $K$ are displayed in Table~\ref{tab:freeFlyingRobotData}.  It is observed from Table~\ref{tab:freeFlyingRobotNum} that both the bicomplex-step, hyper-dual, and automatic differentiation methods consistently require fewer NLP solver iterations to converge than the central finite-difference method using either an over-estimated or exact sparsity pattern.  Furthermore, Table~\ref{tab:freeFlyingRobotTime} shows that both the hyper-dual and  bicomplex-step methods require less computation time than the central finite-difference method because the NLP solver required fewer iterations to converge using either of the former two methods.  Moreover, Table~\ref{tab:freeFlyingRobotItTi} demonstrates how identifying the exact sparsity pattern of the NLP derivative matrices can significantly improve the computational efficiency of the method, as the exact sparsity central finite-difference method takes between $43$ percent to $75$ percent less time on average per NLP solver iteration relative to the over-estimated sparsity central finite-difference method.  Finally, it is observed from Tables~\ref{tab:freeFlyingRobotTime} and \ref{tab:freeFlyingRobotItTi} that automatic differentiation requires substantially more computation time than any of the Taylor series-based methods because automatic differentiation must repeatedly retrace the function in order to obtain the derivative. 

\begin{table}[h!]
\centering
\renewcommand{\baselinestretch}{1}\normalsize\normalfont
\caption{Performance results for Example 1 using $OC$, $EC$, $BC$, $HD$ and $AD$ methods for varying number of mesh intervals $K$ using $N_k=5$ collocation points in each mesh interval. \label{tab:freeFlyingRobotData}}

\subfloat[NLP solver iterations, $I$. \label{tab:freeFlyingRobotNum}]{
\begin{threeparttable}
\begin{tabular}{|c|c|c|c|c|c|}
\hline
$K$	&	$I$&	$\Gamma_{I}^{EC}$	&	$\Gamma_{I}^{BC}$	&	$\Gamma_{I}^{HD}$	&	$\Gamma_{I}^{AD}$	\\ \hline
    2	&	40	&	0	&	22.50	&	22.50	&	22.50	\\
    4	&	45	&	0	&	28.89	&	28.89	&	28.89	\\
    8	&	62	&	16.13	&	19.35	&	19.35	&	19.35	\\
   16	&	48	&	0	&	22.92	&	22.92	&	22.92	\\
   32	&	66	&	0	&	21.21	&	21.21	&	21.21	\\
\hline
\end{tabular}
\end{threeparttable}
}

\subfloat[Computation Time, $\C{T}$. \label{tab:freeFlyingRobotTime}]{
\begin{threeparttable}
\begin{tabular}{|c|c|c|c|c|c|}
\hline
$K$	&	$\C{T}$ (s)	&	$\Gamma_{\C{T}}^{EC}$	&	$\Gamma_{\C{T}}^{BC}$	&	$\Gamma_{\C{T}}^{HD}$	&	$\Gamma_{\C{T}}^{AD}$	\\ \hline
    2	&	0.0876	&	14.96	&	29.54	&	18.30	&	-5982.22	\\
    4	&	0.1205	&	8.54	&	29.28	&	24.48	&	-5741.15	\\
    8	&	0.2510	&	24.51	&	27.44	&	25.50	&	-8390.67	\\
   16	&	0.2973	&	14.43	&	27.34	&	24.76	&	-10326.47	\\
   32	&	0.7769	&	7.32	&	29.36	&	30.71	&	-14885.30	\\
\hline
\end{tabular}
\end{threeparttable}
}

\subfloat[Average Time per NLP Iteration, $\Phi$. \label{tab:freeFlyingRobotItTi}]{
\begin{threeparttable}
\begin{tabular}{|c|c|c|c|c|c|}
\hline
$K$	&	$\Phi$ (ms)	&	$\Gamma_{\Phi}^{EC}$	&	$\Gamma_{\Phi}^{BC}$	&	$\Gamma_{\Phi}^{HD}$	&	$\Gamma_{\Phi}^{AD}$	\\ \hline
    2	&	0.5777	&	64.00	&	68.50	&	64.08	&	-28402.91	\\
    4	&	0.6411	&	59.72	&	61.32	&	56.49	&	-32794.73	\\
    8	&	0.9041	&	57.78	&	51.78	&	52.59	&	-45800.42	\\
   16	&	1.5164	&	48.46	&	41.63	&	38.81	&	-53367.95	\\
   32	&	2.9180	&	35.00	&	30.01	&	34.96	&	-74946.84	\\
\hline
\end{tabular}
\end{threeparttable}
}

\end{table}

%
%
%
%
%
%
%
%

\clearpage

\subsection{Example 2: Minimum Time-to-Climb Supersonic Aircraft Problem \label{subsect:minimumTimeToClimb}}

Consider the following optimal control problem taken from Ref.~\cite{Bryson1}.  Minimize the cost functional  
\begin{equation}\label{eq:minimumTimeToClimbCost}
  \C{J} = t_f~,
\end{equation}
subject to the dynamic constraints
\begin{equation}\label{eq:minimumTimeToClimbDynamics}
\begin{array}{lcl}
  \dot{h}(t) & = & v(t)\sin(\gamma(t))~, \\
  \dot{v}(t) & = & \displaystyle \frac{T(t)\cos(\alpha(t))-D(t)}{m(t)} - \frac{\mu\sin(\gamma(t))}{r^2(t)}~, \\
  \dot{m}(t) & = & \displaystyle -\frac{T(t)}{g_0I_{sp}}~, \\
  \dot{\gamma}(t) & = & \displaystyle \frac{T(t)\sin(\alpha(t))+L(t)}{m(t)v(t)} + \cos(\gamma(t))\left(\frac{v(t)}{r(t)} - \frac{\mu}{v(t)r^2(t)}\right),
\end{array}
\end{equation}
the control inequality constraint
\begin{equation}\label{eq:minimumTimeToClimbControlConstraint}
\begin{array}{c}
\displaystyle -\frac{\pi}{4} \leq \alpha(t) \leq \frac{\pi}{4}~,
\end{array}  
\end{equation}
and the boundary conditions
\begin{equation} \label{eq:minimumTimeToClimbBCs}
\begin{array}{ccccccc}
 h(0) = 0~\text{m}~, & & v(0) = 129.314~\text{m/s}~, & & \gamma(0) = 0~\text{deg}~, & & m(0) = 0~\text{kg}~, \\
 h(t_f) = 19994.88~\text{m}~, & & v(t_f) = 295.092~\text{m/s}~, & & \gamma(t_f) = 0~\text{deg}~, & & m(t_f) = \text{Free}~, \\
\end{array}
\end{equation} 
where 
\begin{equation}
\label{eq:minimumTimeToClimbAux}
\begin{array}{l}
\begin{array}{lclclclclcl}
r(t) & = & h(t) + R_e~,
& &
C_{L_{\alpha}}(t) & = & C_{L_{\alpha}}(M(t))~,
& &
A(t) & = & A(h(t))~,
\\
M(t)~ & = & \displaystyle \frac{v(t)}{A(t)}~,
& &
C_L(t) & = & C_{L_{\alpha}}(t)\alpha(t)~,
& &
T(t) & = & T(h(t),M(t))~,
\\
\rho(t) & = & \rho(h(t))~,
& &
C_{D_0}(t) & = & C_{D_0}(M(t))~,
& &
q(t) & = & \displaystyle \frac{1}{2}\rho(t) v^2(t)~, 
\\
\eta(t) & = & \eta(M(t))~,
&&
L(t) & = & Sq(t)C_L(t)~,
& &
D(t) & = & Sq(t)C_D(t)~,
\end{array}
\\
\begin{array}{lcl}
C_D(t) & = & C_{D_0}(t) + \eta(t)C_{L_{\alpha}}(t)\alpha^2(t)~,
\end{array}
\\
\begin{array}{lclclclclcl}
R_e & = & 6378145~\text{m}~,
& &
\mu & = & 3.986\times10^{14}~\text{m}^3\text{/}\text{s}^{2}~,
\\
g_0~~~~~ & = & 9.80665~\text{m/s}^2~,
&&
S~~ & = & 49.2386~\text{m}^2~,
& &
\hspace{-0.42in}I_{sp} & = & 1600~\text{s}~.
\end{array}
\end{array}
\end{equation}
The models used for the problem were created using data taken from Ref.~\cite{Bryson1}.

The optimal control problem described by Eqs.~\eqref{eq:minimumTimeToClimbCost} -- \eqref{eq:minimumTimeToClimbAux} was approximated using LGR collocation (see Section~\ref{sect:LGR}) using $K=(2,4,8,16,32)$ mesh intervals with $N_k=5$ collocation points in each mesh interval.  The resulting set of NLPs for the different values of $K$ were then solved using an over-estimated and exact sparsity pattern with central finite-difference, and using exact sparsity patterns with bicomplex-step, hyper-dual, and automatic differentiation.  The solutions obtained in all cases are essentially identical and converge to the exact solution (see pages $256$ to $264$ of Ref.~\cite{Betts3}).  Comparisons of the computational efficiency in terms of $I$, $\C{T}$, and $\Phi$ for varying values of $K$ are displayed in Table~\ref{tab:minimumTimeToClimbData}.  It is seen in Table~\ref{tab:minimumTimeToClimbNum} that in nearly all cases the bicomplex-step, hyper-dual, and automatic differentiation methods require fewer NLP solver iterations to converge than the central finite-difference method using either an over-estimated or exact sparsity pattern.  Furthermore, Table~\ref{tab:minimumTimeToClimbTime} shows that both the hyper-dual and bicomplex-step methods require less computation time than the central finite-difference method because the NLP solver requires fewer iterations to converge using either of the former two methods.  Moreover, Table~\ref{tab:minimumTimeToClimbItTi} shows that the bicomplex-step and hyper-dual methods can be more computationally efficient per NLP iteration than the central finite-difference method.  In particular, for this example the central finite-difference method takes approximately $30$ percent longer on average than either of the former two methods.  Finally, it may be observed from Tables~\ref{tab:minimumTimeToClimbTime} and \ref{tab:minimumTimeToClimbItTi} that the $AD$ method requires substantially more computation time than any of the Taylor series-based methods because automatic differentiation must repeatedly retrace the function in order to obtain the derivative. 

\begin{table}[h!]
\centering
\renewcommand{\baselinestretch}{1}\normalsize\normalfont
\caption{Performance results for Example 2 using $OC$, $EC$, $BC$, $HD$, and $AD$ methods for varying number of mesh intervals $K$ using $N_k=5$ in each mesh interval. \label{tab:minimumTimeToClimbData}}

\subfloat[NLP Solver Iterations, $I$. \label{tab:minimumTimeToClimbNum}]{
\begin{threeparttable}
\begin{tabular}{|c|c|c|c|c|c|}
\hline
$K$	&	$I$&	$\Gamma_{I}^{EC}$	&	$\Gamma_{I}^{BC}$	&	$\Gamma_{I}^{HD}$	&	$\Gamma_{I}^{AD}$	\\ \hline
    2	&	248	&	26.61	&	79.03	&	66.94	&	71.37	\\
    4	&	30	&	6.67	&	-43.33	&	-43.33	&	-43.33	\\
    8	&	52	&	-1.92	&	30.77	&	30.77	&	30.77	\\
   16	&	70	&	-7.14	&	32.86	&	32.86	&	32.86	\\
   32	&	53	&	0	&	26.42	&	26.42	&	26.42	\\
\hline
\end{tabular}
\end{threeparttable}
}

\subfloat[Computation Time, $\C{T}$. \label{tab:minimumTimeToClimbTime}]{
\begin{threeparttable}
\begin{tabular}{|c|c|c|c|c|c|}
\hline
$K$	&	$\C{T}$ (s)	&	$\Gamma_{\C{T}}^{EC}$	&	$\Gamma_{\C{T}}^{BC}$	&	$\Gamma_{\C{T}}^{HD}$	&	$\Gamma_{\C{T}}^{AD}$	\\ \hline
    2	&	0.6529	&	31.25	&	81.39	&	69.64	&	-1315.45	\\
    4	&	0.1428	&	5.09	&	-14.69	&	-12.51	&	-7375.47	\\
    8	&	0.4029	&	-1.34	&	44.02	&	41.11	&	-4682.36	\\
   16	&	1.0057	&	-1.12	&	49.06	&	48.54	&	-4541.37	\\
   32	&	1.4072	&	0.71	&	40.26	&	42.88	&	-5130.10	\\
\hline
\end{tabular}
\end{threeparttable}
}

\subfloat[Average Time per NLP Iteration, $\Phi$. \label{tab:minimumTimeToClimbItTi}]{
\begin{threeparttable}
\begin{tabular}{|c|c|c|c|c|c|}
\hline
$K$	&	$\Phi$ (ms)	&	$\Gamma_{\Phi}^{EC}$	&	$\Gamma_{\Phi}^{BC}$	&	$\Gamma_{\Phi}^{HD}$	&	$\Gamma_{\Phi}^{AD}$	\\ \hline
    2	&	1.5267	&	8.21	&	37.41	&	33.18	&	-8227.16	\\
    4	&	2.7067	&	5.92 	&	32.05 	&	36.20	&	-8781.22	\\
    8	&	5.2132	&	2.38 	&	36.38	&	34.10	&	-9851.17	\\
   16	&	10.8110	&	7.92	&	38.30	&	38.23	&	-8865.56	\\
   32	&	20.3869	&	2.69 	&	32.18	&	36.26	&	-8895.19	\\
\hline
\end{tabular}
\end{threeparttable}
}

\end{table}

\clearpage

\subsection{Example 3: Space Station Attitude Control \label{subsect:spaceStation}}

Consider the following optimal control problem taken from Refs.~\cite{Betts3} and \cite{Pietz2003}.  Minimize the cost functional
\begin{equation}\label{eq:space-station-cost}
  \C{J} = {\textstyle\frac{1}{2}}\int_{t_0}^{t_f} \m{u}^{\tr}(t)\m{u}(t)~dt~,
\end{equation}
subject to the dynamic constraints
\begin{equation}\label{eq:space-station-dynamics}
  \begin{array}{lcl}
    \dot{\g{\omega}}(t) & = & \m{J}^{-1}\left\{\g{\tau}_{gg}(\m{r}(t))-\g{\omega}(t)^{\otimes}\left[\m{J}\g{\omega}(t)+\m{h}(t)\right]-\m{u}(t)\right\}~,\vspace{3pt}\\
    \dot{\m{r}}(t) & = &   {\textstyle\frac{1}{2}}\left[\m{r}(t)\m{r}^{\tr}(t)+\m{I}+\m{r}(t)\right]\left[\g{\omega}(t)-\g{\omega}_0(\m{r}(t))\right]~,\vspace{3pt}\\
    \dot{\m{h}}(t) & = & \m{u}(t)~,
  \end{array}
\end{equation}
the inequality path constraint
\begin{equation}
  \left\|\m{h}(t)\right\| \leq h_{\max}~,
\end{equation}
and the boundary conditions
\begin{equation}\label{eqLspace-station-boundary-conditions}
\begin{array}{c}
  \begin{array}{rclcrcl} 
    t_0 & = & 0~, & &
    t_f & = & 1800~,\quad\quad\quad\quad\quad\quad\quad\quad\quad\quad\quad~ \\
    \g{\omega}(0) & = & \bar{\g{\omega}}_0~, & &
    \m{r}(0) & = & \bar{\m{r}}_0~, \vspace{3pt}\\
    \m{h}(0) & = & \bar{\m{h}}_0~, & &
  \end{array} \\
  \begin{array}{rcl}
    \mathbf{0} & = & \mathbf{J}^{-1}\left\{\boldsymbol{\tau}_{gg}(\mathbf{r}(t_f)) - \boldsymbol{\omega}^{\otimes}(t_f)\left[\mathbf{J}\boldsymbol{\omega}(t_f)+\mathbf{h}(t_f)\right]\right\}~, \vspace{3pt}\\
     \mathbf{0} & = & {\textstyle\frac{1}{2}}\left[\mathbf{r}(t_f)\mathbf{r}^{\textsf{T}}(t_f)+\mathbf{I}+\mathbf{r}(t_f)\right]\left[\boldsymbol{\omega}(t_f)-\g{\omega}_0(\mathbf{r}(t_f))\right]~,
  \end{array}
\end{array}
\end{equation}
where $(\g{\omega}(t),\m{r}(t),\m{h}(t))$ is the state, $\m{u}(t)$ is the control, $\m{I}$ is the identity matrix, and $h_{max}=10000$.  In this formulation, $\g{\omega}(t)$ is the angular velocity, $\m{r}(t)$ is the Euler-Rodrigues parameter vector, $\m{h}(t)$ is the angular momentum, and $\m{u}(t)$ is the input moment (and is the control).  Additionally, 
\begin{equation}\label{eq:space-station-omega-orb-defs}
  \begin{array}{rclcrcl}
    \g{\omega}_0(\m{r}(t)) & = & -\omega_{\textrm{orb}}\m{C}_2(\m{r}(t))~, & &
    \g{\tau}_{gg}(\m{r}(t)) & = & 3\omega_{\textrm{orb}}^2 \m{C}_3^{\otimes}(\m{r}(t))\m{J}\m{C}_3(\m{r}(t))~,
  \end{array}
\end{equation}
and $\m{C}_2(\m{r}(t))$ and $\m{C}_3(\m{r}(t))$ are the second and third column, respectively, of the matrix
\begin{equation}\label{eq:space-station-C-def}
 \m{C}(\m{r}(t)) = \m{l}+\frac{2}{1+\m{r}^{\tr}(t)\m{r}(t)}\left(\m{r}^{\otimes}(t)\m{r}^{\otimes}(t)-\m{r}^{\otimes}(t)\right)~,
\end{equation}
where $\omega_{orb}=0.06511\pi/180$.  In this example the matrix $\m{J}$ is given as
\begin{equation}\label{eq:space-station-inertia-matrix}
  \m{J} = \left[
    \begin{array}{ccc} 2.80701911616\times 10^{7} & 4.822509936\times 10^{5} & -1.71675094448\times 10^{7}\\
      4.822509936\times 10^{5} & 9.5144639344\times 10^{7}  & 6.02604448\times 10^{4} \\
      -1.71675094448\times 10^{7} & 6.02604448\times 10^{4} & 7.6594401336\times 10^{7}
    \end{array}
  \right]~,
\end{equation}
while the initial conditions $\bar{\g{\omega}}_0$, $\bar{\m{r}}_0$, and
$\bar{\m{h}}_0$ are, respectively,
\begin{equation}\label{eq:space-station-initial-conditions}
  \begin{array}{lll}
    \bar{\g{\omega}}_0 = \left[\begin{array}{c} -9.53807\times 10^{-6}\\ -1.13633\times 10^{-3} \\+5.34728\times 10^{-6}\end{array}\right]~,\quad
\bar{\m{r}}_0 = \left[\begin{array}{c} 2.99637\times10^{-3}  \\1.53345\times 10^{-1} \\ 3.83598\times 10^{-3}\end{array}\right]~,~\text{and}~
\bar{\m{h}}_0 = \left[\begin{array}{c} 5000 \\ 5000 \\ 5000 \end{array}\right]~.
  \end{array}
\end{equation}
A more detailed description of this problem, including all of the constants $\m{J}$, $\bar{\g{\omega}}_0$, $\bar{\m{r}}_0$, and $\bar{\m{h}}_0$, can be found in either Ref.~\cite{Betts3} or Ref.~\cite{Pietz2003}.  

The optimal control problem described in Eqs.~\eqref{eq:space-station-cost} -- \eqref{eq:space-station-initial-conditions} was approximated using LGR collocation (see Section~\ref{sect:LGR}) using $K=(2,4,8,16,32)$ mesh intervals with $N_k=5$ collocation points in each mesh interval.  The resulting set of NLPs for the different values of $K$ were then solved using an over-estimated and exact sparsity pattern with central finite-difference, and an exact sparsity pattern using bicomplex-step, hyper-dual, and automatic differentiation.  The solutions obtained in all cases are essentially identical and converge to the exact solution (see pages $296$ -- $297$ of Ref.~\cite{Betts3}).  Comparisons of the computational efficiency in terms of $I$, $\C{T}$, and $\Phi$ for varying values of $K$ are displayed in Table~\ref{tab:spaceStationData}.   It is seen in Table~\ref{tab:spaceStationNum} that for $K<32$ mesh intervals that the number of NLP solver iterations is approximately the same for any of the five methods.  On the other hand, for $K=32$ mesh intervals, both the bicomplex-step and hyper-dual methods require $39.72$ percent fewer iterations than the over-estimated central finite-difference method, while the exact sparsity central finite-difference and automatic differentiation methods require $24.11$ percent and $36.88$ percent fewer iterations, respectively.  Furthermore, Table~\ref{tab:spaceStationTime} shows that for all of the cases studied both the hyper-dual and bicomplex-step methods also require less computation time than the over-estimated sparsity central finite-difference method.  Moreover, Table~\ref{tab:spaceStationItTi} shows the increase in computational efficiency per iteration when using an exact NLP sparsity pattern, as the exact sparsity central finite-difference method takes $61$ percent to $73$ percent less time on average per NLP iteration than the over-estimated sparsity central finite-difference method on all tested instances.  Finally, it is observed in Tables~\ref{tab:spaceStationTime} and \ref{tab:spaceStationItTi} that the automatic differentiation method requires substantially more computation time than any of the Taylor series-based methods.  

\begin{table}[h!]
\centering
\renewcommand{\baselinestretch}{1}\normalsize\normalfont
\caption{Performance results for Example 3 using $OC$, $EC$, $BC$, $HD$, and $AD$ methods for varying number of mesh intervals $K$ with $N_k=5$ collocation points in each mesh interval.\label{tab:spaceStationData}}

\subfloat[NLP Solver Iterations, $I$. \label{tab:spaceStationNum}]{
\begin{threeparttable}
\begin{tabular}{|c|c|c|c|c|c|}
\hline
$K$	&	$I$&	$\Gamma_{I}^{EC}$	&	$\Gamma_{I}^{BC}$	&	$\Gamma_{I}^{HD}$	&	$\Gamma_{I}^{AD}$	\\ \hline
    2	&	22	&	0	&	0	&	0	&	0	\\
    4	&	33	&	0	&	0	&	0	&	0	\\
    8	&	44	&	-4.55	&	-4.55	&	-4.55	&	-4.55	\\
   16	&	50	&	-6.00	&	-10.00	&	-6.00	&	-8.00	\\
   32	&	141	&	24.11		&	39.72	&	39.72	&	36.88	\\
\hline
\end{tabular}
\end{threeparttable}
}

\subfloat[Computation Time, $\C{T}$. \label{tab:spaceStationTime}]{
\begin{threeparttable}
\begin{tabular}{|c|c|c|c|c|c|}
\hline
$K$	&	$\C{T}$ (s)	&	$\Gamma_{\C{T}}^{EC}$	&	$\Gamma_{\C{T}}^{BC}$	&	$\Gamma_{\C{T}}^{HD}$	&	$\Gamma_{\C{T}}^{AD}$	\\ \hline
    2	&	0.2628	&	53.71	&	55.53	&	52.39 	&	-7420.02	\\
    4	&	0.4703	&	50.86	&	46.98	&	49.26	&	-11200.24	\\
    8	&	0.8522	&	43.72	&	33.75 	&	33.63	&	-16350.86	\\
   16	&	1.4118	&	40.61	&	20.10	&	27.41	&	-22562.69	\\
   32	&	7.7361	&	53.20	&	57.38	&	57.50	&	-12922.78	\\
\hline
\end{tabular}
\end{threeparttable}
}

\subfloat[Average Time per NLP Iteration, $\Phi$. \label{tab:spaceStationItTi}]{
\begin{threeparttable}
\begin{tabular}{|c|c|c|c|c|c|}
\hline
$K$	&	$\Phi$ (ms)	&	$\Gamma_{\Phi}^{EC}$	&	$\Gamma_{\Phi}^{BC}$	&	$\Gamma_{\Phi}^{HD}$	&	$\Gamma_{\Phi}^{AD}$	\\ \hline
    2	&	8.2873	&	73.38	&	74.45	&	75.30	&	-10222.58	\\
    4	&	9.2951	&	72.13	&	68.57	&	69.24	&	-16673.23	\\
    8	&	11.6082	&	69.75	&	54.49	&	56.43	&	-25547.80	\\
   16	&	15.8056	&	65.05	&	36.13	&	40.89	&	-36645.81	\\
   32	&	26.1159	&	60.06	&	26.06	&	26.66	&	-42707.56	\\
\hline
\end{tabular}
\end{threeparttable}
}

\end{table}

%
%
%
%
%
%
%
%


\section{Discussion \label{sect:Discussion}}

As developed in Sections~\ref{sect:multi-phase} through \ref{sect:NLP}, transcribing a continuous optimal control problem into a large sparse nonlinear programming problem (NLP) using direct collocation introduces the necessity of an efficient method for the NLP solver employed to solve the resulting NLP.  The central finite-differencing, bicomplex-step derivative approximation, and hyper-dual derivative approximation described in Sections~\ref{sect:finite-difference} through \ref{sect:HyperDual} have truncation error estimates based on Taylor series expansions \cite{Gill3,Lantoine1,Fike2011}.  Section~\ref{subsect:Deriv-Approx-Error} demonstrated that, while the central finite-differencing and the bicomplex-step derivative approximations have truncation errors on the order of $\mathcal{O}(h^2)$, the bicomplex-step derivative approximation has a significant advantage over central finite-differencing because of the roundoff error associated with central finite differencing.  Moreover, as seen in Fig.~\ref{fig:ExampleFuncDerErr}, the bicomplex-step derivative approximation can be employed using an arbitrarily small step size such that the truncation error reaches near machine precision \cite{Lantoine1}.  Conversely, when employing central finite-differencing, a step size must be chosen in a manner that minimizes the sum of the truncation error and the roundoff error \cite{Gill3}.  Next, the hyper-dual derivative approximation does not suffer from either truncation error or roundoff error.  As a result, an arbitrary step size may be used when implementing the hyper-dual derivative approximation, and this approximation is always accurate to machine precision.  The accuracy of the hyper-dual derivative approximation is also shown in Fig.~\ref{fig:ExampleFuncDerErr}.  Finally, the derivative estimates obtained using the automatic differentiation discussed in Section~\ref{sect:AlgorithmicDiff} are accurate to machine precision, because the derivative is computed algorithmically using the calculus chain rule.  

The immunity of the bicomplex-step and hyper-dual derivative approximations to roundoff error enables either of these derivative approximation methods to provide extremely accurate derivative approximations for the NLP solver.  In fact, the derivative estimates obtained using either of these methods are comparable in accuracy to the derivatives obtained using automatic differentiation.  The accuracy of the derivative estimates supplied to the NLP solver has an enormous impact on the search direction taken by the NLP solver and the associated rate of convergence.  As demonstrated by the three example problems shown in Section~\ref{sect:Examples}, the increased accuracy of the derivative estimates provided by the bicomplex-step, hyper-dual, and automatic differentiation methods most often enables the NLP solver to converge in fewer iterations when compared with using the derivative estimates obtained using central finite-differencing.  The reason that the NLP solver does not in all cases converge in fewer iterations using bicomplex-step, hyper-dual, or automatic differentiation when compared with central finite-differencing is due to the fact that the NLP solver employs a Newton method to determine the search direction on each iteration.  In particular, the Newton method is based upon a quadratic approximation which may not lead to the fastest convergence rate in the case where the NLP is highly nonlinear.  Thus even though the bicomplex-step, hyper-dual, and automatic differentiation methods provide more accurate approximations of the NLP derivative matrices, the search directions chosen upon each iteration using the less accurate central finite-difference methods may, in a relatively infrequent number of cases, inadvertently cause the NLP solver to converge in fewer iterations.

For all three example problems in Section~\ref{sect:Examples}, the fewer number of NLP solver iterations leads to reduced computation time when using either the bicomplex-step or hyper-dual method as compared to the central finite-difference method.  This lower computation time arises from the fact that the majority of computation time to solve the NLP is spent computing derivative approximations.  These fewer iterations lead to fewer computations of the derivative estimates and makes the bicomplex-step and hyper-dual methods advantageous over the central finite-difference method.  Interestingly, the fewer number of iterations required when using automatic differentiation is outweighed by the large computation time per iteration required when using automatic differentiation relative to the other derivative approximation methods.

Finally, for all three Taylor series-based derivative approximations presented in Sections~\ref{sect:finite-difference} through \ref{sect:HyperDual}, a perturbation to the inputs of the function whose derivative is of interest is required in order to compute the derivative approximation.  As discussed in Section~\ref{sect:finite-difference}, central finite-differencing is subject to both roundoff and truncation error, thus limiting the minimum step size that can be used.  The effect of the less accurate derivative approximation of central finite-differencing is demonstrated in Section~\ref{sect:Examples} where the central finite-difference methods generally require more NLP solver iterations to converge when compared with either bicomplex-step, hyper-dual, or automatic differentiation.  While the bicomplex-step derivative approximation described in Section~\ref{sect:Bicomplex} appears to be subject to only truncation error, in practice roundoff error can occur due to the evaluation of the bicomplex number arithmetic when using double precision computations.  The fact that some of the derivative approximations become less accurate for step sizes below a certain magnitude when implementing the bicomplex-step derivative approximation is an artifact of using double precision arithmetic in the implementation.  On the other hand, the hyper-dual derivative approximation presented in Section~\ref{sect:HyperDual} is immune to both roundoff and truncation error.  Thus when implementing the hyper-dual derivative approximation, the perturbation step size is arbitrary and does not affect the accuracy of the approximation.  Conversely, the accuracy of central finite-differencing and bicomplex-step derivative approximations are dependent upon the step size, leading to scaling and rounding issues that make implementation more difficult and less accurate when compared with the hyper-dual approximation.  Finally, although automatic differentiation provides derivative estimates that are accurate to machine precision, the computational overhead required is found to be excessively expensive, thus making implementation impractical when compared to any of the Taylor series-based derivative approximation methods.

\section{Conclusions \label{sect:Conclusions}}

Four derivative estimation methods, central finite-differencing, bicomplex-step, hyper-dual, and automatic differentiation, have been compared in terms of their effectiveness for use with direct collocation methods for solving optimal control problems.  The process of transcribing a continuous optimal control problem into a large sparse nonlinear programming problem using a previously developed Legendre-Gauss-Radau direct collocation method is described.  The form of the resultant nonlinear programming problem and the need for an efficient method to facilitate the NLP solver employed is presented.  The three Taylor series-based derivative approximations, central finite-differencing, bicomplex-step, and hyper-dual, are derived.  These three Taylor series-based methods are then compared to one another and compared to automatic differentiation in terms of accuracy, computational efficiency, and implementation.  The performance of a nonlinear programming problem solver is then demonstrated on three benchmark optimal control problems.  The performance of the NLP solver is assessed in terms of number of iterations to solve, computation time per iteration, and computation time.  Despite the observation that central finite-differencing requires less computation time per iteration than either the bicomplex-step or the hyper-dual method, the latter two methods require significantly less overall computation time because the NLP solver requires significantly fewer iterations to converge to a solution.  Moreover, the bicomplex-step and hyper-dual methods are found to have similar performance, although the hyper-dual method is found to be significantly easier to implement.  Additionally, automatic differentiation is found to require substantially more computation time than any of the Taylor series-based methods.  Finally, a preliminary comparison of the derivative estimation methods for solving optimal control problems using direct collocation is found to favor the hyper-dual method in terms of computational efficiency and ease of implementation.

\section*{Acknowledgments}

The authors gratefully acknowledge support for this research from the U.S.~Office of Naval Research under grant N00014-15-1-2048 and from the U.S.~National Science Foundation under grants DMS-1522629, DMS-1924762, and CMMI-1563225.  

\renewcommand{\baselinestretch}{1.0}
\normalsize\normalfont
\bibliographystyle{aiaa}

\begin{thebibliography}{10}
\newcommand{\enquote}[1]{``#1''}

\bibitem{Betts3}
Betts, J.~T., {\em Practical Methods for Optimal Control and Estimation Using
  Nonlinear Programming\/}, SIAM Press, Philadelphia, 2nd ed., 2009.

\bibitem{Gill3}
Gill, P.~E., Murray, W., and Wright, M.~H., {\em Practical Optimization\/},
  Academic Press, London, 1981.

\bibitem{Biegler2}
Biegler, L.~T. and Zavala, V.~M., \enquote{{L}arge-{S}cale {N}onlinear
  {P}rogramming {U}sing {IPOPT}: {A}n {I}ntegrating {F}ramework for
  {E}nterprise-{W}ide {O}ptimization,} {\em Computers and Chemical
  Engineering\/}, Vol.~33, No.~3, March 2008, pp.~575--582.

\bibitem{Gill1}
Gill, P.~E., Murray, W., and Saunders, M.~A., \enquote{{SNOPT}: {A}n {SQP}
  {A}lgorithm for {L}arge-{S}cale {C}onstrained {O}ptimization,} {\em SIAM
  Review\/}, Vol.~47, No.~1, January 2002, pp.~99--131.
  \url{https://doi.org/10.1137/S0036144504446096}.

\bibitem{Martins2}
Martins, J.~R. and Hwang, J.~T., \enquote{Review and Unification of Methods for
  Computing Derivatives of Multidisciplinary Computational Models,} {\em AIAA
  Journal\/}, Vol.~51, No.~11, September 2013, pp.~2582--2599.
  \url{https://doi.org/10.2514/1.J052184}.

\bibitem{Sobieszczanski1}
Sobieszczanski-Sobieski, J., \enquote{Sensitivity of Complex, Internally
  Coupled Systems,} {\em AIAA Journal\/}, Vol.~28, No.~1, January 1990,
  pp.~153--160. \url{https://doi.org/10.2514/3.10366}.

\bibitem{Griewank0}
Griewank, A. and Walther, A., {\em Evaluating Derivatives: Principles and
  Techniques of Algorithmic Differentiation\/}, SIAM Press, Philadelphia,
  Pennsylvania, 2nd ed., 2008.

\bibitem{Griewank1}
Griewank, A., Juedes, D., and Utke, J., \enquote{{Algorithm 755}: {ADOL-C}: A
  Package for the Automatic Differentiation of Algorithms Written in {C\slash
  C++},} {\em {ACM} Transactions on Mathematical Software\/}, Vol.~22, No.~2,
  1996, pp.~131--167. \url{https://doi.org/10.1145/229473.229474}.

\bibitem{Walther1}
Walther, A., Griewank, A., and Vogel, O., \enquote{ADOL-C: Automatic
  Differentiation Using Operator Overloading in C++,} {\em Proceedings in
  Applied Mathematics and Mechanics\/}, Vol.~2, No.~1, 2003, pp.~41--44.
  \url{https://doi.org/10.1002/pamm.200310011}.

\bibitem{Betts13}
Betts, J.~T. and Huffman, W.~P., \enquote{{E}xploiting {S}parsity in the
  {D}irect {T}ranscription {M}ethod for {O}ptimal {C}ontrol,} {\em
  Computational Optimization and Applications\/}, Vol.~14, 1999, pp.~179--201.
  \url{https://doi.org/10.1023/A:1008739131724}.

\bibitem{Patterson2012}
Patterson, M.~A. and Rao, A.~V., \enquote{{E}xploiting {S}parsity in {D}irect
  {C}ollocation {P}seudospectral {M}ethods for {S}olving {C}ontinuous-{T}ime
  {O}ptimal {C}ontrol {P}roblems,} {\em Journal of Spacecraft and Rockets,\/},
  Vol.~49, No.~2, March--April 2012, pp.~364--377.
  \url{https://doi.org/10.2514/1.A32071}.

\bibitem{Agamawi2018}
Agamawi, Y.~M. and Rao, A.~V., \enquote{{E}xploiting {S}parsity in {D}irect
  {O}rthogonal {C}ollocation {M}ethods for {S}olving {M}ultiple-{P}hase
  {O}ptimal {C}ontrol {P}roblems,} {\em 2018 Space Flight Mechanics Meeting\/},
  AIAA Paper2018-0724, Kissimmee, Florida, 8-12 January.
  \url{https://doi.org/10.2514/6.2018-0724} 2018.

\bibitem{Betts2}
Betts, J.~T., \enquote{Survey of Numerical Methods for Trajectory
  Optimization,} {\em Journal of Guidnance, Control, and Dynamics\/}, Vol.~21,
  No.~2, March--April 1998, pp.~193--207. \url{https://doi.org/10.2514/2.4231}.

\bibitem{Lantoine1}
Lantoine, G., Russell, R.~P., and Dargent, T., \enquote{{U}sing {M}ulticomplex
  {V}ariables for {A}utomatic {C}omputation of {H}igh-{O}rder {D}erivatives,}
  {\em ACM Transactions on Mathematical Software\/}, Vol.~38, April 2012,
  pp.~16:1--16:21. \url{10.1145/2168773.2168774}.

\bibitem{Fike2011}
Fike, J. and Alonso, J., \enquote{The Development of Hyper-Dual Numbers for
  Exact Second-Derivative Calculations,} {\em 49th AIAA Aerospace Sciences
  Meeting including the New Horizons Forum and Aerospace Exposition\/}, AIAA
  Paper 2011-886, Orlando, Florida, 4-7 January 2011.
  \url{https://doi.org/10.2514/6.2011-886}.

\bibitem{Garg1}
Garg, D., Patterson, M.~A., Darby, C.~L., Francolin, C., Huntington, G.~T.,
  Hager, W.~W., and Rao, A.~V., \enquote{Direct Trajectory Optimization and
  Costate Estimation of Finite-Horizon and Infinite-Horizon Optimal Control
  Problems via a Radau Pseudospectral Method,} {\em Computational Optimization
  and Applications\/}, Vol.~49, No.~2, June 2011, pp.~335--358.
  \url{https://doi.org/10.1007/s10589--00--09291--0}.

\bibitem{Garg2}
Garg, D., Patterson, M.~A., Hager, W.~W., Rao, A.~V., Benson, D.~A., and
  Huntington, G.~T., \enquote{A Unified Framework for the Numerical Solution of
  Optimal Control Problems Using Pseudospectral Methods,} {\em Automatica\/},
  Vol.~46, No.~11, November 2010, pp.~1843--1851.
  \url{https://doi.org/10.1016/j.automatica.2010.06.048}.

\bibitem{Garg3}
Garg, D., Hager, W.~W., and Rao, A.~V., \enquote{Pseudospectral Methods for
  Solving Infinite-Horizon Optimal Control Problems,} {\em Automatica\/},
  Vol.~47, No.~4, April 2011, pp.~829--837.
  \url{https://doi.org/10.1016/j.automatica.2011.01.085}.

\bibitem{Kameswaran1}
Kameswaran, S. and Biegler, L.~T., \enquote{Convergence Rates for Direct
  Transcription of Optimal Control Problems Using Collocation at Radau Points,}
  {\em Computational Optimization and Applications\/}, Vol.~41, No.~1, 2008,
  pp.~81--126. \url{https://doi.org/10.1007/s10589--007--9098--9}.

\bibitem{Patterson2015}
Patterson, M.~A., Hager, W.~W., and Rao, A.~V., \enquote{A $ph$ Mesh Refinement
  Method for Optimal Control,} {\em Optimal Control Applications and
  Methods\/}, Vol.~36, No.~4, July--August 2015, pp.~398--421.
  \url{https://doi.org/10.1002/oca.2114}.

\bibitem{Liu2015}
Liu, F., Hager, W.~W., and Rao, A.~V., \enquote{{A}daptive {M}esh {R}efinement
  for {O}ptimal {C}ontrol {U}sing {N}onsmoothness {D}etection and {M}esh {S}ize
  {R}eduction,} {\em Journal of the Franklin Institute\/}, Vol.~352, No.~10,
  October 2015, pp.~4081--4106.
  \url{https://doi.org/10.1016/j.jfranklin.2015.05.028}.

\bibitem{Liu2018}
Liu, F., Hager, W.~W., and Rao, A.~V., \enquote{{A}daptive {M}esh {R}efinement
  for {O}ptimal {C}ontrol {U}sing {U}sing {D}ecay {R}ates of {L}egendre
  {P}olynomial {C}oefficients,} {\em IEEE Transactions on Control System
  Technology\/}, Vol.~26, No.~4, 2018, pp.~1475 -- 1483.
  \url{https://doi.org/10.1109/TCST.2017.2702122}.

\bibitem{HagerHouRao15a}
Hager, W.~W., Hou, H., and Rao, A.~V., \enquote{Lebesgue Constants Arising in a
  Class of Collocation Methods,} {\em IMA Journal of Numerical Analysis\/},
  Vol.~13, No.~1, October 2017, pp.~1884--1901.
  \url{https://doi.org/10.1093/imanum/drw060}.

\bibitem{HagerHouRao16a}
Hager, W.~W., Hou, H., and Rao, A.~V., \enquote{Convergence Rate for a Gauss
  Collocation Method Applied to Unconstrained Optimal Control,} {\em Journal of
  Optimization Theory and Applications\/}, Vol.~169, No.~3, 2016, pp.~801--824.
  \url{10.1007/s10957--016--0929--7}.

\bibitem{HagerLiuMohapatraWangRao18a}
Hager, W.~W., Liu, J., Mohapatra, S., Rao, A.~V., and Wang, X.-S.,
  \enquote{Convergence Rate for a Gauss Collocation Method Applied to
  Constrained Optimal Control,} {\em SIAM Journal on Control and
  Optimization\/}, Vol.~56, No.~2, 2018, pp.~1386 -- 1411.
  \url{https://doi.org/10.1137/16M1096761}.

\bibitem{Abramowitz1}
Abramowitz, M. and Stegun, I., {\em Handbook of Mathematical Functions with
  Formulas, Graphs, and Mathematical Tables\/}, Dover Publications, New York,
  1965.

\bibitem{Patterson2014}
Patterson, M.~A. and Rao, A.~V., \enquote{$\mathbb{GPOPS-II}$, {A} {MATLAB}
  {S}oftware for {S}olving {M}ultiple-{P}hase {O}ptimal {C}ontrol {P}roblems
  {U}sing $hp$-{A}daptive {G}aussian {Q}uadrature {C}ollocation {M}ethods and
  {S}parse {N}onlinear {P}rogramming,} {\em ACM Transactions on Mathematical
  Software\/}, Vol.~41, No.~1, October 2014, pp.~1:1--1:37.
  \url{https://doi.org/10.1145/2558904}.

\bibitem{Sakawa}
Sakawa, Y., \enquote{Trajectory Planning of a Free-Flying Robot by Using the
  Optimal Control,} {\em Optimal Control Applications and Methods\/}, Vol.~20,
  1999, pp.~235--248.
  \url{https://doi.org/10.1002/(SICI)1099--1514(199909/10)20:5<235::AID--OCA658>3.0.CO;2--I}.

\bibitem{Bryson1}
Bryson, A.~E., Desai, M.~N., and Hoffman, W.~C., \enquote{Energy-State
  Approximation in Performance Optimization of Supersonic Aircraft,} {\em AIAA
  Journal of Aircraft\/}, Vol.~6, No.~6, November--December 1969, pp.~481--488.
  \url{https://doi.org/10.2514/3.44093}.

\bibitem{Pietz2003}
Pietz, J.~A., {\em Pseudospectral Collocation Methods for the Direct
  Transcription of Optimal Control Problems\/}, Master's thesis, Rice
  University, Houston, Texas, April 2003.

\end{thebibliography}

\end{document}